\documentclass[twoside,a4paper,leqno]{siamltex}
\usepackage{amsmath
,amssymb,fullpage}
\usepackage{graphicx}
\usepackage{hyperref}
\hypersetup{%
  pdftitle={},
  pdfsubject={},
  pdfauthor={},
  pdfkeywords={},
  pdfcreator={pdfLaTeX},
  pdfproducer={pdfLaTeX},
  }

\usepackage{color}
\usepackage{subfigure}

\newcommand{\N}{\mathbf{N}}
\newcommand{\R}{\mathbf{R}}
\newcommand{\Z}{\mathbf{Z}}

\newtheorem{remark}[equation]{Remark}
\begin{document}
\title{Pulsatile Viscous Flows in Elliptical Vessels and Annuli:
  Solution to the Inverse Problem, with Application to Blood and
  Cerebrospinal Fluid Flow}
\author{Luigi C.\ Berselli\thanks{Dipartimento di Matematica,
    Universit{\`a} di Pisa, Via F.~Buonarroti 1/c, I-56127 Pisa,
    ITALY.  ({berselli@dma.unipi.it})} \and Francesca
  Guerra\footnotemark[2] \and Barbara Mazzolai,\footnotemark[2] \and
  Edoardo Sinibaldi\footnotemark[2] }

\date{\today}\maketitle
\renewcommand{\thefootnote}{\fnsymbol{footnote}}
\footnotetext[2]{Center for Micro-BioRobotics@SSSA, Istituto Italiano di
Tecnologia, Viale R. Piaggio 34, Pontedera, Italy, (edoardo.sinibaldi@iit.it)}
\begin{abstract}
  We consider the fully-developed flow of an incompressible Newtonian
  fluid in a cylindrical vessel with elliptical cross-section (both an
  ellipse and the annulus between two confocal ellipses).  In
  particular, we address an \textit{inverse problem}, namely to
  compute the velocity field associated with a given, time-periodic
  flow rate.  This is motivated by the fact that flow rate is the main
  physical quantity which can be actually measured in many practical
  situations.  We propose a novel numerical strategy, which is
  nonetheless grounded on several analytical relations.  The proposed
  method leads to the solution of some simple ordinary differential
  systems. It holds promise to be more amenable to implementation than
  previous approaches, which are substantially based on the
  challenging computation of Mathieu functions.  Some numerical
  results are reported, based on measured data for human blood flow in
  the internal carotid artery, and cerebrospinal fluid (CSF) flow in
  the upper cervical region of the human spine.  As expected,
  computational efficiency is the main asset of our solution: a
  speed-up factor over $10^3$ was obtained, compared to more elaborate
  numerical approaches.

  The main goal of the present study is to provide an improved source
  of initial/boundary data for more ambitious numerical approaches, as
  well as a benchmark solution for pulsatile flows in elliptical
  sections with given flow rate.  The proposed method can be
  effectively applied to bio-fluid dynamics investigations (possibly
  addressing key aspects of relevant diseases), to biomedical
  applications (including targeted drug delivery and energy harvesting
  for implantable devices), up to longer-term medical microrobotics
  applications.
\end{abstract}
%
\section{Introduction}
\label{sec_intro}
Pulsatile flows, here intended as flows driven by a time-periodic
force (generally the pressure gradient), are encountered in many
contexts of practical interest.  Amongst them, a remarkable case is
that of heart-beat driven, human physiological flows, including blood
circulation~\cite{rogers2010blood} and (less directly) cerebrospinal
fluid (CSF) flow~\cite{Irani2008,wagshul2011the}.  Besides bio-fluid
dynamics, time-periodic flows are widely studied with regard to
e.g. chemical-physics applications, like species
separation~\cite{Huang2006}, as well as for mass and heat transfer
problems~\cite{Korichi2009,Wang2005}.  Furthermore, many studies can
be found in the literature which address pulsatile flows in a more
general (yet still application-driven) fluidic context, like
peristaltic pumping~\cite{kumar2010}, to name but one.  However, in
many cases of practical interest the pressure gradient is unknown or
hardly measurable, and the only available datum is the time-dependent
flux (hereafter, flux and flow rate are understood as synonyms).  This
is the most common measurement, for instance, in clinical settings
regarding blood and CSF flow.

As for blood flow, it is difficult for many portions of the
vasculature to keep an ideal, circular cross-section, due to the
presence of surrounding organs (possibly activated by relevant
muscles, or simply displaced by postural changes).  Hence, while
keeping some degree of abstractness, it is more appropriate to
consider an elliptical cross-section to obtain a more accurate
description of the flow (see e.g.~\cite{HZ1998}). As for CSF flow, the
spinal sub-arachnoid space can be well approximated by an elliptical
annulus~\cite{Irani2008}; CSF dynamics in such a domain is affected by
pulsatility (mainly towards the upper cervical region) and plays a
major role in the (still poorly understood) pathophysiology of many
high-impact disease, like syringomyelia and Chiari
malformation~\cite{dirocco2011hyd,loukas2011}.  Of course, in order to
tackle realistic 3D flow conditions, for instance associated with
patient-specific geometries, a fully numerical approach is mandatory
(see
e.g.~\cite{cheng2012the,Hentschel2010,Linn2009,LYA2001,Shaffer2011}
for the CSF case). This implies the exploitation of rather
time-consuming approaches, even when adopting a classical Newtonian
behavior for the working fluid (as surely appropriate for
CSF~\cite{LYA2001}, as well as for blood in large
vessels~\cite{FQV2006,RSO2009}).  Nevertheless, it is possible to keep
some degree of physical representativeness also by adopting the
hypothesis of a fully-developed flow (see e.g.~\cite{GPK2008}),
despite its inherent oversimplifications from a physical
viewpoint. Such a position has been systematically adopted in many
studies in order to obtain analytical flow solutions (see below),
otherwise hardly achievable.  However, even by considering a
fully-developed flow, it is necessary to solve a non-standard,
\textit{inverse problem} in order to evaluate the velocity field (and
the pressure gradient) by starting from a given flux across any
cross-section of a vessel.

With reference to the above observations, we point out to expressly
confine our attention to a special class of flow solutions, namely
fully-developed flows of an incompressible Newtonian fluid in a
straight cylinder with elliptical cross-section (either simply
connected or not, see below). In particular, these assumptions permit
to address a simplified linear problem, for which superposition
principle holds.  Nevertheless, the adopted positions aim at obtaining
a benchmark solution for the inverse problem described below: this is
the main target of the present work. More precisely, by virtue of the
fully-developed flow hypothesis, we provide a numerical strategy which
is substantially grounded on analytical relations, while being more
amenable to implementation than previous ones proposed in
literature.  Hence, while laying no claims of generality (not even
deliberately addressing realistic three-dimensional flows, e.g. within
patient specific geometries and possibly involving complex rheological
effects), the method we propose is credited to hold some value for
obtaining at least an approximation to real flow conditions, while
containing the computational cost, in the spirit of benchmark flow
solutions.  Furthermore, the obtained solution strategy can be used as
an improved source of boundary data for more ambitious numerical
approaches based on realistic data, up to serving as a debugging tool.
It is worth stressing once again that we address the inverse problem,
for which no straightforward approaches have been well established so
far.

To better define the problem, we consider the motion of an
incompressible (the constant density is here normalized to 1),
Newtonian fluid in a semi-infinite straight pipe $P=E\times
\R^+\subset\R^3$. In the sequel $E\subset \R^2$ will be either an
ellipse or an elliptical annulus (even if results for circular
cross-sections are also recalled, and some abstract results can be
obtained for any smooth and bounded cross-section).  The governing
equations are therefore the incompressible Navier-Stokes equations
which, in a reference frame with $z$ directed along the pipe axis (and
$x:=(x_1,x_2)$ belonging to an orthogonal plane), read
\begin{equation*}
  \begin{aligned}
    \partial_t \vec{u}+(\vec{u}\cdot\nabla)\,\vec{u}-\nu\Delta
    \vec{u}+\nabla p=0&\qquad (x,z)\in E\times\R^+,\; t\in\R^+,
    \\
    \nabla\cdot \vec{u}=0&\qquad (x,z)\in E\times\R^+,\; t\in\R^+,
    \\
    \vec{u}=0&\qquad (x,z)\in \partial E\times\R^+,\; t\in\R^+,
  \end{aligned}
\end{equation*}
where $\vec{u}(t,x,z)$ and $p(t,x,z)$ respectively denote velocity and
pressure, and $\nu$ represents kinematic viscosity.  Furthermore, as
anticipated, we deliberately neglect turbulent effects by looking for
a special class of solutions, namely fully-developed solutions (also
named Poiseuille-type solutions) of the following form:
\begin{equation*}
  \vec{u}(t,x,z) = (0,0,w(t,x))
    \qquad\text{and}\qquad
  p(t,x,z) = -\lambda(t,x,z)+p_0(t),
\end{equation*}
where $p_0(t)$ denotes an arbitrary function which disappears from the equations,
since pressure only enters through its spatial gradient in the
formulation of the problem. In addition, the following flux condition is
assumed to be satisfied:
\begin{equation*}
  \int_{E} w(t,x)\,dx=f(t),
\end{equation*}
for some given scalar function $f(t)$, only depending on the time
variable.
By standard arguments, the above Poiseuille-type \textit{ansatz} implies that
%
%
we may assume that the pressure has the form
$p(t,z)=-\lambda(t)\, z$. Moreover, the dependence of $w$ on the space
variables $x_1$ and $x_2$ allows us to consider a problem reduced to
the cross-section $E$, with the given flux condition.  This
implies that we have to study the following problem: Find $(w(t,x),\lambda(t))$ such
that
\begin{equation}
  \label{eq:basic_flow}
  \begin{cases}
    \partial_t w(t,x) - \nu\Delta_x w(t,x) = \lambda(t),&\qquad x\in
    E,\; t\in\R^+,
    \\
    w(t,x)=0&\qquad x\in \partial E,\; t\in\R^+,
    \\
    \int_{E} w(t,x)\,dx=f(t)&\qquad t\in\R^+,
  \end{cases}
\end{equation}
where $\Delta_x$ denotes the Laplacian with respect to the variables
$x_1$ and $x_2$.  This can be considered as an inverse problem since,
if the axial pressure gradient $\lambda$ is known, one can evaluate
$w$ (by solving a two-dimensional scalar heat equation) and thus the
corresponding flux as well. Contrarily, in our problem the datum is
the flux $f(t)$ and one wants to find the corresponding velocity and
pressure associated with the fully-developed solution (such a problem
is linked to one of the nowadays classical Leray's
problems~\cite{Gal1994a}).  For the sake of completeness, we recall
that existence and uniqueness of the fully-developed solution in pipes
with very general cross-sections (and with explicit but not exact
expression in terms of the flux) have been recently obtained
in~\cite{Bei2005c,BR2010,Pil2007}, under reasonable regularity on the
periodic $f(t)$.

The relevance of the cross-section shape on the solution process
deserves some discussion at this stage.  Indeed, when considering a
circular cross-section, one can write an explicit solution in terms of
special functions. In particular, the analytical solution of the
direct problem was firstly obtained by Sexl~\cite{Sex1930} and
Womersley~\cite{Wom1955} in terms of Bessel functions, while we
recently derived in~\cite{BS2011} an analytical solution to the
inverse problem involving regularized confluent hyper-geometric
functions.  Such a fully analytical approach is made possible by the
particularly simple expression taken by the Laplace operator in
cylindrical coordinates. In fact, when considering an elliptical
cross-section, in spite of the apparent similarity with the circular
case, the situation drastically changes: an explicit solution can only
be obtained for the stationary problem, while it is necessary to
resort to numerical computations already when addressing the direct
time-dependent
problem~\cite{HZ1998,Kha1957,RD2004,Verma1960ellAnn,Verma1960ell}.  In
particular, current solutions in literature are based on involved
tools like the Mathieu functions~\cite{McLach1947}.  Nevertheless,
some numerical approaches which also tackle the inverse problem have
been recently proposed~\cite{GPK2008}, still by invoking the Mathieu
functions, also involving the determination of the Laplacian
eigenvalues in an elliptical domain.  Yet this is not a trivial task:
there are several works emphasizing how instabilities can arise in the
numerical computations when addressing such a problem by directly
managing the Mathieu functions, especially when the ellipse
eccentricity is small~\cite{HZ1998}.  Indeed, computing the Mathieu
functions has been recognized as the most critical issue
in~\cite{GPK2008}; major efforts have been spent in order to define a
wise algorithm for their evaluation.

In the light of these points, we propose an alternative numerical
approach to the inverse problem of a pulsatile fully-developed flow in
elliptical vessels.  Our approach only indirectly addresses the
Mathieu problem and therefore it does not suffer from the
aforementioned limitations. Some numerical simulations are also
reported, based on flow rates actually measured for human blood flow
in the internal carotid artery, and CSF flow in the upper cervical
region of the human spine.  Corresponding results support the
effective usability of our formulation (successfully coping with small
eccentricities as well), by virtue of the dramatic reduction in
computational time it allows, as compared to more elaborated finite
element approaches. This, in turn, broadens the application spectrum
of our solution, as highlighted in the last section of the present
work.

\smallskip

\textbf{Plan of the paper:} In Section~\ref{sec_prob_staz} we review
relevant issues regarding stationary, fully-developed flows in
elliptical cross-sections, either simply connected or not.  Some
results related to circular cross-sections are also recalled, for
completeness. In Section~\ref{sec:time-dependent} we propose a novel,
computationally-efficient numerical strategy for solving the inverse
time-dependent problem, by also exploiting an auxiliary, direct
formulation.  In particular, besides recalling recent results for the
the circular cross-section, we discuss in detail the case of the
simply connected elliptical section (being the annular case very
similar). Furthermore, exemplificative numerical simulations are
reported in Section~\ref{sec_numtest}, based on blood and CSF
physiological flow data; computational cost of the proposed strategy
is also compared to more elaborate, finite element numerical
approaches.  Concluding remarks are finally reported in
Section~\ref{sec_conclusions}, where the potential for effective
application of our solution is highlighted, with explicit references
to the biomedical and microrobotic research field.

\section{Stationary problem}
\label{sec_prob_staz}

Despite being classical, the stationary solution in the elliptical
cross-section is re-derived in a way that is functional to subsequent
treatment.  Classical results for the circular cross-section are
recalled as well, for completeness.  Such a stationary problem is
formulated as a Poisson problem with Dirichlet boundary conditions
(besides the flux condition), namely: Find $(w(x),\lambda)$ such that
\begin{equation}
  \label{eq:basic_flow-stationary}
  \begin{cases}
  -\nu\Delta_x w(x) = \lambda,&\qquad x\in
    E,
    \\
    w(x)=0,&\qquad x\in \partial E,
    \\
    \int_{E} w(x)\,dx=f,&
  \end{cases}
\end{equation}
where $f\in\R$ is the given flux, while $\left(-\lambda\right)\in\R$ is the constant
pressure gradient.
%
\subsection{Stationary flow in a circular cross-section}
The history of this problem dates back to the work of Hagen and
Poiseuille (see e.g.~\cite{Gal1994a} for details), who firstly derived
the solution in a circular domain.  Nonetheless, the interest for such
a special solution (which can be easily derived thanks to complete
separation of variables) is still considerable, as a benchmark flow
also in applications to turbulence (see e.g.~\cite{Pop2000}), or for
innovative biomedical applications (like magnetic particle
targeting~\cite{Haverkort2009}).  This is mainly due to the fact that
Poiseuille flow represents one of the few examples of exact solution
of fluid equations which are not space-periodic, yet with Dirichlet
boundary conditions.

\subsubsection{Flow in the circle}
Let $E$ be a circle
with radius $R>0$ and centered at the origin, namely:
$E=\left\{(x_1,x_2)\in\R^2:\ |x|<R\right\}$,
where $|x|=\sqrt{x_1^2+x_2^2}$.
For such a case the solution to~\eqref{eq:basic_flow-stationary} is
easily obtained, namely the well-known parabolic profile
\begin{equation}
  \label{eq:Poiseuille-circular}
  w(x)=\frac{2 f}{\pi R^2}\left(1-\frac{|x|^2}{R^2}\right).
\end{equation}
%
\subsubsection{Flow in the circular annulus}
Let $E$ be a circular annulus delimited by radii $R_1$ and $R_2>R_1$, namely:
$E:=\left\{(x_1,x_2)\in\R^2:\ R_1<|x|<R_2\right\}$.
In this case the solution of the Poisson problem~\eqref{eq:basic_flow-stationary} is
\begin{equation}
w(x)=  \frac{2 f}{\pi \, R_2^2}
\frac{\left(1 - \frac{|x|^2}{R_2^2}\right) -
\left(1 - \frac{R_1^2}{R_2^2}\right)
\frac{\log\left(\frac{|x|}{R_2}\right)}{\log\left(\frac{R_1}{R_2}\right)}}
{\left(1 - \frac{R_1^4}{R_2^4}\right) +
\left(1 - \frac{R_1^2}{R_2^2}\right)^2
\frac{1}{\log\left(\frac{R_1}{R_2}\right)}}.
\label{sol_staz_ann_circ}
\end{equation}
Notice the limiting behavior of~\eqref{sol_staz_ann_circ} for $R_2=R$ and
$R_1\to0$.

\subsection{Stationary flow in an elliptical cross-section}
We review some well-known results, which will be needed later on to
handle the time-dependent case. Relevant solutions for this problem
originally appeared in the pioneering works of Khamrui~\cite{Kha1957}
and Verma~\cite{Verma1960ell,Verma1960ellAnn} (besides earlier,
preliminary contributions cited therein).  Also in the present case it
is possible to obtain an analytical solution, thanks to the
particularly simple expression of the external force in the Poisson
problem (\ref{eq:basic_flow-stationary}). Nonetheless, it is worth
mentioning that numerical approximation of the aforementioned problem
in an elliptical domain, yet with more general external forces, is the
subject of ongoing research (see e.g.~\cite{Lai2004}).

\subsubsection{Flow in the ellipse}
\label{lab_flusso_ell_mathieau}
Let $E=\left\{(x_1,x_2)\in\R^2:
  x^2_1/\alpha^2+x^2_2/\beta^2<1\right\}$ denote an ellipse, where
$\alpha=a\cosh(b)$ and $\beta=a\sinh(b)$ respectively denote the
length of the major ($x_1$-direction) and minor semi-axis
($x_2$-direction), while $2a$ represents the inter-focal distance.  We
preliminarily observe that the equality
\begin{equation}
  \cosh^{-1}\left(\frac{\alpha}{a}\right)=b=\sinh^{-1}\left(\frac{\beta}{a}\right)
\label{eq_b_fromalfabeta} 
\end{equation}
implies $\log\left(\alpha/a+\sqrt{\alpha^2-a^2}/a\right)=
\log\left(\beta/a +\sqrt{\beta^2 +a^2}/a\right)$, from which
$b=\log\left(\left(\alpha+\beta\right)/\sqrt{\alpha^2-\beta^2}\right)$
and $a=\sqrt{\alpha^2-\beta^2}$.  Then, to construct the solution we
introduce the natural change of variables
\begin{equation}
  \label{eq:change_variable}
  \begin{cases}
    &    x_1=a \cosh(\eta)\cos(\theta),
    \\
    &  x_2=a \sinh(\eta)\sin(\theta),
  \end{cases}
\end{equation}
with Jacobian
\begin{equation}
\label{eq_J_etatheta}
  J(\eta,\theta) = a^2 \big[\sinh^2(\eta)+\sin^2(\theta)\big]
  = \frac{a^2}{2} \left[-\frac{e^{-2i\theta}}{2}+
                         \cosh (2\eta )
                        -\frac{e^{2i\theta}}{2}\right].
\end{equation}
By employing this change of variables we can reduce the problem in the
rectangular domain
\begin{equation}
  E'=\left\{(\eta,\theta)\in \R^2:
  \ 0<\eta<b\quad\text{and}\quad 0\leq\theta<2\pi\right\},
\label{eq_rett_dom_ellipt}
\end{equation}
and it is well-known that the Laplace operator is separable with these
coordinates.  Furthermore, by applying the aforementioned change of
variables, the Poisson equation~\eqref{eq:basic_flow-stationary}
becomes
\begin{equation*}
  \left\{
  \begin{aligned}
    -\frac{\nu}{
      J(\eta,\theta)}\Bigg(\frac{\partial^2 u}{\partial \eta^2}+
    \frac{\partial^2 u}{\partial
      \theta^2}\Bigg)&=\lambda\qquad(\eta,\theta)\in E',
    \\
    u\text{ is $2\pi$-periodic in $\theta$}&
    \\
    u(b,\theta)&=0\qquad\theta\in[0,2\pi[,
    \\
  \int_{E'} u(\eta,\theta)\,J(\eta,\theta)\,d\eta \,d\theta=f
  \end{aligned}
\right.
\end{equation*}
where $u(\eta,\theta)$ simply represents $w(x)$ in the new variables
(please also notice that $J>0$ for $\eta>0$).  Moreover, we observe
that an additional Neumann condition is needed at $\eta=0$, as
detailed below. Due to periodicity in $\theta$, the solution can be
written with the following complex Fourier series expansion with
respect to $\theta$ (hereafter ``hat'' $\hat{\cdot}$ denotes a complex
Fourier coefficient):
\begin{equation*}
  u(\eta,\theta)=\sum_{n\in\Z} \widehat{u}_n(\eta)\,e^{in\theta},
  \qquad \mathrm{where} \qquad 
  \widehat{u}_n(\eta)=\frac{1}{2\pi}\int_0^{2\pi} u(\eta,\theta)\,e^{-in\theta}.
\end{equation*}
By virtue of linearity, we temporarily assume $\lambda=1$ and
we recast the problem at hand as follows (for the function $u$ which is
$2\pi$-periodic in $\theta$ and temporarily we do not assign the flux):
\begin{equation*}
  \left\{
    \begin{aligned}
      - {\nu}\Big(\frac{\partial^2 u}{\partial \eta^2}+
      \frac{\partial^2 u}{\partial
        \theta^2}\Big)&=J(\eta,\theta)
      \qquad&(\eta,\theta)\in E',
      \\
      u(b,\theta)&=0\qquad &\theta\in[0,2\pi[.
    \end{aligned}
  \right.
\end{equation*}
We then plug in the Poisson equation the Fourier series for both $u$
and $J(\eta,\theta)$, as given in~\eqref{eq_J_etatheta}, and we
observe that, in order to have a real-valued solution, we need to
impose as usual $\widehat{u}_{-n}=\overline{\widehat{ u}_{n}}$
(hereafter ``bar'' $\bar{\cdot}$ denotes complex conjugacy).  In
particular, $\Im(\widehat{u}_0)=0$ (hereafter $\Re$ and $\Im$
respectively denote the real and imaginary part).  Moreover, by using
the symmetries of the solution (reflections over the ellipse axes of
symmetry) it follows that $\widehat{u}_n=0$ for $n$ odd.  Furthermore,
in order to have a smooth solution also at $\eta=0$, we need to impose
$\partial_\eta u(0,\theta)=0$.  As a result, by exploiting the above
conditions, and by observing that $J$ only possesses three non-zero
modes (namely those associated with $n=-2,0,2$), we obtain the
following identity:
\begin{equation*}
  -\nu\sum_{n=-2,0,2}  [\widehat{u}_n''(\eta)-n^2\widehat{u}_n(\eta)]
  \,e^{in\theta}=\frac{a^2}{2} \left[-\frac{e^{-2i\theta}}{2}+ 
    \cosh (2\eta )-\frac{e^{2i\theta}}{2}\right],
\end{equation*}
which corresponds to a system of three complex ordinary differential
equations. In the above expression and in the sequel, ``prime''
$({\,\cdot\,}')$ denotes derivative with respect to $\eta$.  Then, by
equating the terms corresponding to the same exponential, and by
recalling the conditions at $\eta=0,b$ coming from the Dirichlet
condition and from the symmetries of the solution, we have to solve
the following boundary value problems:

\begin{equation*}
  \left\{  \begin{aligned}
      &\widehat{u}_{\pm2}''(\eta)-4\widehat{u}_{\pm2}
      =\frac{a^2}{4\nu},
      \\
      &\widehat{u}_{\pm2}(b)=\Re(\widehat{u}_{\pm2}'(0))=\Im(\widehat{u}_{\pm2}(0))= 0,
    \end{aligned}
  \right.\qquad\text{and}\qquad
  \left\{
    \begin{aligned}
      &\widehat{u}_0''(\eta) =-\frac{a^2}{2\nu} \cosh (2 \eta ) ,
      \\
      &\widehat{u}_{0}(b)=\widehat{u}_{0}'(0)=0 .
  \end{aligned}
\right.
\end{equation*}
Hence, by explicitly solving these uncoupled linear differential
equations we get
\begin{equation*}
  \widehat{u}_0(\eta) = -\frac{a^2}{8 \nu} \left(\cosh(2 \eta)-\cosh(2 b)\right),
  \qquad
  \widehat{u}_{\pm2}(\eta)=-\frac{a^2}{16 \nu}
                            \frac{1+e^{4 b}-e^{2 b-2 \eta}-e^{2 b+2\eta}}{1+e^{4 b}}. 
\end{equation*}
Then, by summing-up the solutions after multiplication by the
corresponding complex exponential, and after some algebraic
manipulations, the solution to the problem with $\lambda=1$ turns out
to be
\begin{equation*}
  u(\eta,\theta)=\frac{a^2}{8 \nu }
                \left( \cosh (2 b)-\cos (2 \theta) \right) \,
                \left( \cosh (2 b)-\cosh (2 \eta)  \right) \,
                \text{sech}(2 b), 
\end{equation*}
being the corresponding flux
\begin{equation*}
  \int_{E'} u(\eta,\theta)\,J(\eta,\theta)\,d\eta \,d\theta
  = a^2 \, \int_{0}^{b}
     \left(\int_{0}^{2\pi} u(\eta,\theta) \, [\sinh^2(\eta)+\sin^2(\theta)] \, d\theta \right) \, d\eta 
     = \frac{\pi a^4}{32 \nu} \, \sinh^2(2b) \, \tanh(2b).
\end{equation*}
That said, the solution to the original
problem~\eqref{eq:basic_flow-stationary} with given flux $f$ is
finally obtained as $\left(\lambda u(\eta,\theta),\lambda\right)$,
i.e. by scaling the function $u(\eta,\theta)$ computed right above by
the following factor:
\begin{equation}
  \lambda=f \, \frac{{32 \nu}}{\pi a^4} \left[\sinh^2(2 b) \, \tanh(2 b)\right]^{-1},
\label{eq_f2lambda_staz}
\end{equation}
accounting for the actual flow rate.  Finally, by mapping back to Cartesian
coordinates, the sought solution reads
\begin{equation}
  w(x_1,x_2)=\frac{2 f}{\pi  \alpha  \beta}
  \left[1-\frac{x^2_1}{\alpha^2}-\frac{x_2^2}{\beta^2}\right].
\label{staz_ellisse_scalato}
\end{equation}
It is interesting to note how the
solution~\eqref{staz_ellisse_scalato} derives from Poiseuille
solution~\eqref{eq:Poiseuille-circular} by a simple anisotropic
scaling of the variables (as for obtaining the considered elliptical
domain by starting from the unit circle).  Thus, the elegant
expression~\eqref{staz_ellisse_scalato} (appearing in
e.g.~\cite{HZ1998}, yet missing in many related works) could have been
obtained through simpler derivations. Nevertheless, we deliberately
decided to introduce many details above, since they serve as useful
references for the solution of the corresponding time-dependent
problem.  On regard, it should be noticed how, despite separation of
variables (in elliptical coordinates), ellipticity \textit{per se}
creates a new situation already when considering stationary
conditions: there appear 3 coupled modes in the solution (contrarily
to the fully decoupled flow solution in the circle).  Such a
difference persists and it is somehow magnified in the time-dependent
case, where more complex behaviors naturally appear (see
Section~\ref{sec:time-dependent}).

\subsubsection{Flow in the elliptical annulus}
\label{subsub_ellip_annul}
Given the inter-focal distance $2a$, it is straightforward to
apply~\eqref{eq:change_variable} in order to also map two confocal
ellipses into a rectangular region.  More precisely, given semi-axes
$\alpha_1,\beta_1$ and $\alpha_2,\beta_2$, with $\beta_2>\beta_1$, the
rectangle $E'=\left\{(\eta,\theta)\in \R^2: \
  b_1<\eta<b_2\quad\text{and}\quad 0\leq\theta<2\pi\right\}$ is
obtained, where $b_1$ and $b_2 > b_1$ are derived from the respective
semi-axes through~\eqref{eq_b_fromalfabeta}.  For this problem we have
Dirichlet boundary conditions at both $\eta=b_1$ and
$\eta=b_2$. Hence, by enforcing the symmetries
$\widehat{u}_{-n}=\overline{\widehat{u}_n}$ needed to have real-valued
velocities, we end up with the following boundary value problems (for
the auxiliary problem with $\lambda=1$):
\begin{equation*}
  \left\{  \begin{aligned}
      &\widehat{u}_{\pm2}''(\eta)-4\widehat{u}_{\pm2}
      =\frac{a^2}{4\nu},
      \\
      &\widehat{u}_{\pm2}(b_1)=\widehat{u}_{\pm2}(b_2)= 0,
    \end{aligned}
  \right.\qquad\text{and}\qquad
  \left\{
    \begin{aligned}
      &\widehat{u}_0''(\eta) =-\frac{a^2}{2\nu} \cosh (2 \eta )
      \\
      &\widehat{u}_{0}(b_1)=\widehat{u}_{0}(b_2)=0,
  \end{aligned}
\right.
\end{equation*}
whose solutions, by direct integration, read:
\begin{eqnarray*}
 \widehat{u}_0(\eta) & = & \frac{a^2}{8 \nu} 
\left[\frac{(\eta-b_1)\cosh(2b_2)-(\eta-b_2)\cosh(2b_1)}{(b_2-b_1)}-\cosh(2 \eta)\right],\\
\widehat{u}_{\pm2}(\eta) & = & \frac{a^2}{16 \nu}
\left[ \frac{e^{2(\eta-b_2)} - e^{-2(\eta-b_2)} - e^{2(\eta-b_1)} + e^{-2(\eta-b_1)} }{e^{2(b_1-b_2)} - e^{-2(b_1-b_2)}}
-1 \right], 
\end{eqnarray*}
so that the following stationary solution for the auxiliary problem is
obtained (cf.~\cite{Verma1960ellAnn}):
\begin{eqnarray*}
  u(\eta,\theta) & = & \frac{a^2}{8 \nu}\left[
    \frac{(\eta-b_1)\cosh(2b_2)-(\eta-b_2)\cosh(2b_1)}{(b_2-b_1)}\; + \right. \\
  & & \qquad \quad \left. 
    \frac{\sinh(2(\eta-b_1)) - \sinh(2(\eta-b_2))}{\sinh(b_2-b_1)}\cos(2\theta)  -(\cosh(2\eta)+\cos(2\theta))\right]. 
\end{eqnarray*}
Finally, by reasoning as in Section~\ref{lab_flusso_ell_mathieau},
the solution to~\eqref{eq:basic_flow-stationary}
with given flux $f$ is obtained as $\left(\lambda u(\eta,\theta),\lambda\right)$,
with (cf.~\cite{GPK2008})
\begin{equation*}
 \lambda = f \, \frac{16 \nu}{\pi a^4}
\left[  \frac{(\sinh(4 b_2)-\sinh(4 b_1))}{4}
      - \frac{(\cosh(2 b_2)-\cosh(2 b_1))^2}{2(b_2-b_1)}
      - \frac{\cosh(2(b_2-b_1))-1}{\sinh(2(b_2-b_1))} \right]^{-1}.
\end{equation*}
Following~\cite{Verma1960ellAnn}, it is worth remarking that the
solution in the (simply connected) ellipse cannot be obtained as a
limit case of the one in the annulus, due to the constraint of
confocality. Indeed, by allowing $b_1\to0$, the inner boundary of the
annulus (where the Dirichlet condition also applies) degenerates to
the inter-focal segment. The resulting flow field is therefore that
one associated with an ellipse with semi-axes $\alpha_2$ and
$\beta_2$, yet also containing a plate having width equal to the
inter-focal distance $2a$.

\section{Time-dependent problem}
\label{sec:time-dependent}
We now address the time-periodic elliptical case, by proposing a novel
numerical method for its solution (this is the main asset of the
present study).  More in detail, we firstly recall the solution of the
inverse problem for the flow in the circle, so as to keep some degree
of symmetry with Section~\ref{sec_prob_staz}.  We then address the
elliptical case, where we solve the inverse problem by exploiting
relevant results, purposely introduced through an auxiliary, direct
formulation.  In particular, we accurately report the solution method
for the simply connected elliptical cross-section. Conversely, we omit
details for the flow in the elliptical annulus between two confocal
ellipses, since it can be derived from the previous case by a
straightforward modification of the boundary conditions, as mentioned
in Section~\ref{subsub_ellip_annul}; corresponding numerical results are
nonetheless reported in Section~\ref{sec_res_CSF}, for the sake of
completeness.
%
\subsection{Flow in the circle: the inverse problem}
Let us preliminarily introduce the following $T$-periodic functions:
\begin{equation}
  \lambda(t)=\sum_{m\in\Z}\widehat{\lambda}_m e^{i\omega_m t}, \quad
        f(t)=\sum_{m\in\Z}\widehat{f}_m e^{i\omega_m t}, \quad
        w(x,t)=\sum_{m\in\Z}\widehat{w}_m(x) \,e^{i \omega_m t}, \quad\quad
      \omega_m=\frac{2\pi m}{T}.
\label{def_omegam}
\end{equation}
The aforementioned works by Sexl~\cite{Sex1930} and
Womersley~\cite{Wom1955} addressed the direct problem, that is with
assigned pressure gradient $\lambda(t)$.  Corresponding solution can
be constructed as an infinite series of Bessel functions $J_0$ of
order zero (Womersley considered a single Fourier mode, yet
superposition can be invoked under the assumptions in~\cite{Bei2005c},
since the corresponding sequence formally converges).  As regards the
inverse problem, an explicit map between the Fourier coefficients of
the velocity $\widehat{w}_m$ and those of the flow rate
$\widehat{f}_m$ has been recently obtained by the
authors~\cite{BS2011}.  Relevant relations read
\begin{equation*}
  \widehat{f}_m =  \pi R^2  \,
  \left(1-\frac{\,_0\tilde{F}_1\left(;2;i \, \mathrm{Wo}^2_{R,m} \, /4 \right)}
               {\,_0\tilde{F}_1\left(;1;i \, \mathrm{Wo}^2_{R,m} \, /4 \right)}
  \right) \,\, \frac{\widehat{\lambda}_m}{i \omega_m},
  \qquad
  \widehat{w}_m= \left(1-\frac{J_0\left( (-1)^{3/4} \, \mathrm{Wo}_{r,m} \right)}
                              {J_0\left( (-1)^{3/4} \, \mathrm{Wo}_{R,m} \right)}
                 \right) \,\, \frac{\widehat{\lambda}_m}{i \omega_m},       
\end{equation*}
where $\,_0\tilde{F}_1(;\cdot;\cdot)$ denotes the
regularized confluent hyper-geometric function, and
\begin{equation} 
\label{wo_num}
\mathrm{Wo}_{r,m}= r \, \sqrt{\frac{\omega_m}{\nu}},
\end{equation}
is a non-dimensional parameter which generalizes the classical
Womersley number.
%
%
\subsection{Flow in the ellipse: the auxiliary direct problem}
\label{sez_probdir_aux}

The explicit solution for the direct problem in an elliptical vessel
was originally derived by Khamrui~\cite{Kha1957} (although some
details were missing) and Verma~\cite{Verma1960ell}, following the
same approach of Womersley.  In particular, given a single harmonic
pressure gradient $\lambda= e^{i\omega_m t}$, with $m\in\Z$, the
solution (in elliptical coordinates) is the following:

\begin{equation}
  u(\eta,\theta)=\sum_{n=0}^\infty C_{2n} \,\, {C\!e}_{2n}(\eta,-q)
  \, {ce}_{2n}(\theta,-q),
  \quad \mathrm{with} \quad q=\frac{i \, a^2 \, \omega_m}{4 \, \nu}.
\label{eq_sol_diretto_conq}
\end{equation}

\noindent
The functions $C\!e_{2n}$ and ${ce}_{2n}$
in~\eqref{eq_sol_diretto_conq} respectively denote the ordinary and
modified Mathieu functions~\cite{McLach1947}, while $C_{2n}$ represent
suitable constants, determined from the no-slip (i.e. Dirichlet)
boundary condition.  It must be remarked that, despite the
aforementioned functions were originally introduced by Mathieu in the
19th century (to study the vibration of an elliptical membrane, and
they are still applied in the field~\cite{CMZ1994}), they still
deserve attention from a computational viewpoint, since their
evaluation is prone to numerical instabilities~\cite{SW2009}.
Moreover, in order to get a $\theta$-periodic solution through the
Mathieu functions (as in the present case), one needs to evaluate the
eigenvalues of the Laplacian: this involves infinite tridiagonal
matrices and it is very expensive.  Furthermore, practical limitations
on the ellipticity $\epsilon=\beta/\alpha$ arise when directly using
these special functions; see e.g.~\cite{HZ1998}, where numerical
instabilities are reported for $\epsilon<0.3$ (besides interesting
expressions for $\epsilon \to 1$, that is for elliptical
cross-sections degenerating to the circular geometry).  In light of
these points, the quest for robust numerical methods able to
efficiently compute the Mathieu functions is still an open issue, and
positive contributions are being provided.  In particular, a stepwise
procedure is proposed in~\cite{GPK2008}, suitable for evaluating the
Mathieu functions in correspondence to large complex arguments
typically associated with viscous flows. Such a strategy exploits a
proper blend of backward and forward recurrence techniques in order to
enhance the convergence of the involved computations.

In order to circumvent the aforementioned limitations, we propose
a novel numerical strategy, in the spirit of Fourier analysis, which
only involves the Mathieu functions in an \textit{indirect way}. In
particular, we extend the approach introduced in~\cite{Lai2004} for
the stationary case to the time-dependent one: our approach is based
on a Fourier analysis in the variables $\theta$ and $t$, while the
problem in the variable $\eta$ is kept in the physical space.  More in
detail, by applying the change of variables introduced
in~\eqref{eq:change_variable}, we firstly recast the time-dependent
problem as follows:
\begin{equation}
  u_t(t,\eta,\theta) \, J(\eta,\theta)
  -\nu \left(\frac{\partial^2 u(t,\eta,\theta)}{\partial \eta^2}+
             \frac{\partial^2 u(t,\eta,\theta)}{\partial \theta^2}\right) = 
  J(\eta,\theta) \, \lambda(t).
\label{eq_auxdirprob_ge1}
\end{equation}
Then, since we look for $T$-periodic solutions (besides needing
$2\pi$-periodicity with respect to $\theta$), we make the following
\textit{ansatz} of separation of variables for the velocity:
\begin{equation}
  u(t,\eta,\theta)=\sum_{m,n\in\Z}\widehat{u}_{m,n}(\eta)\,e^{in\theta}
  e^{i\omega_m t},
\label{eq_ansatz_doppiof}
\end{equation}
hereafter considered in place of the corresponding one
in~\eqref{def_omegam}.  Furthermore, since we consider real-valued
solutions, the following conditions immediately apply:
$\widehat{\lambda}_{-m}=\overline{\widehat{\lambda}_{m}}$,
$\widehat{f}_{-m}=\overline{\widehat{f}_m}$ and
$\widehat{u}_{-m,-n}=\overline{\widehat{u}_{m,n}}$ (latter relation
implies, in particular, $\Im(\widehat{u}(\eta)_{0,0})=0$).

We then start by deliberately addressing a direct problem, hence by
assuming a \textit{given} pressure gradient.  In order to stress the
fact that we temporarily turn our attention to a direct problem, we
introduce a minor notation change (just within this section) by
denoting the pressure gradient as
$\sigma(t)=\sum_{m\in\Z}\widehat{\sigma}_m e^{i\omega_m t}$.  By
plugging~\eqref{eq_ansatz_doppiof} into~\eqref{eq_auxdirprob_ge1}, we
arrive at the following problem: Solve, for each $m\in \Z$,
\begin{equation}
  \left\{
    \begin{aligned}
      \sum_{n\in\Z} i\,\omega_m \widehat{u}_{m,n}(\eta) \,e^{in\theta}
      J(\eta,\theta)
      -{\nu}\big(\widehat{u}_{m,n}''(\eta)-n^2
      \widehat{u}_{m,n}(\eta)\big)\,e^{in\theta}&=J(\eta,\theta)\,\widehat{\sigma}_m
      \\
      \widehat{u}_{m,n}(b)&=0\qquad n\in\Z,
  \end{aligned}
\right.
\label{eq_auxdirprob_ge2}
\end{equation}
with $(\eta,\theta)\in E'$ and $E'$ defined as
in~\eqref{eq_rett_dom_ellipt}.  The considered problem can be
immediately simplified by invoking symmetry: velocity must be
unchanged by the transformations $\theta\mapsto-\theta$ and
$\theta\mapsto\pi-\theta$, i.e. by reflection over the ellipse axes of
symmetry.  This, by also recalling the above conditions on conjugacy,
provides
\begin{equation}
  \widehat{u}_{m,n}=0 \quad \text{for $n$ odd} \quad  \text{and} \quad
  \widehat{u}_{m,-n}=\widehat{u}_{m,n}, \,\,
  \widehat{u}_{-m,n}=\overline{\widehat{u}_{m,-n}}=\overline{\widehat{u}_{m,n}}, \quad
    \text{for $n$ even, $m\in\Z$.}
\label{eq:symmetries}
\end{equation}
Hence, we only need to solve the equations for the Fourier modes
$\widehat{u}_{m,n}$, with $n,m\in\Z$ both greater or equal to zero,
with $n$ even. This permits to reduce by a factor eight the
computational burden.  That said, by plugging
into~\eqref{eq_auxdirprob_ge2} the Fourier expansion of $J$ given
in~\eqref{eq_J_etatheta} and by equating the corresponding even
$n$-modes (at fixed $m$), we obtain the following infinite family of
ordinary differential equations
%
    \begin{equation}
    \begin{array}{rcl}
      \widehat{u}''_{m,0} -
         \bigg[     \frac{ i \, \mathrm{Wo}^2_{a,m}}{2} \cosh(2\eta)\bigg]\,\widehat{u}_{m,0}
                   +\frac{ i \, \mathrm{Wo}^2_{a,m}}{2} \,
              \widehat{u}_{m,2}
         &=&-\frac{a^2}{2 \nu}\cosh(2\eta) \, \widehat{\sigma}_m
      \\
      \widehat{u}''_{m,2} -
         \bigg[2^2 +\frac{ i \, \mathrm{Wo}^2_{a,m}}{2} \cosh(2\eta)\bigg]\,\widehat{u}_{m,2}
                   +\frac{ i \, \mathrm{Wo}^2_{a,m}}{4} \,(\widehat{u}_{m,4} +\widehat{u}_{m,0})
         &=&\frac{a^2}{4 \nu} \, \widehat{\sigma}_m
      \\
      \widehat{u}''_{m,4} -
         \bigg[4^2 +\frac{ i \, \mathrm{Wo}^2_{a,m}}{2} \cosh(2\eta)\bigg]\,\widehat{u}_{m,4}
                   +\frac{ i \, \mathrm{Wo}^2_{a,m}}{4} \,(\widehat{u}_{m,6} +\widehat{u}_{m,2})
         &=&\ \ 0
      \\
      \vdots\hspace{2cm}&\vdots&\quad\vdots
      \\
      \widehat{u}''_{m,2n} -
         \bigg[(2n)^2 +\frac{ i \, \mathrm{Wo}^2_{a,m}}{2} \cosh(2\eta)\bigg]\,\widehat{u}_{m,2n}
                      +\frac{ i \, \mathrm{Wo}^2_{a,m}}{4} \,(\widehat{u}_{m,2n+2}+\widehat{u}_{m,2n-2})
         &=&\ \ 0
      \\
      \vdots\hspace{2cm}&\vdots&\quad\vdots 
    \end{array}
  \label{eq_auxdirprob_ge3}
  \end{equation}
  where we use $\widehat{u}_{m,n}$ in place of
  $\widehat{u}_{m,n}(\eta)$ for brevity, and $\mathrm{Wo}_{a,m}$ is a
  generalized Womersley number based on the semi-focal distance $a$,
  defined as in~(\ref{wo_num}).  Please also notice that $q=i
  \mathrm{Wo}^2_{a,m}/4$, with $q$ introduced
  in~\eqref{eq_sol_diretto_conq} (cf. e.g.~\cite{GPK2008,HZ1998}).
  Furthermore, as in the stationary case, we need to assign an
  additional condition at $\eta=0$ (where the change of coordinates
  becomes singular) in order to have a well-posed mixed boundary value
  problems for $\widehat{u}_{m,n}$. In particular, for the solution to
  be smooth also near $\eta=0$, it is suitable to enforce symmetry for
  the velocity-profile, i.e. a vanishing velocity gradient normal to
  the inter-focal line. This reads (cf.~\cite[Eq.~(9)]{HZ1998})
  $\partial_{\eta}u(0,\theta,t)=0$, that is $\widehat{u}'(0)_{m,n}=0$.

\medskip

The equations~\eqref{eq_auxdirprob_ge3} represent an infinite system
of non-homogeneous Mathieu ordinary differential equations.  In order
to get a computable system, one main idea is to approximate such a
system with a finite dimensional (large enough) system of coupled
ordinary differential equations.  However, this cannot be achieved by
simply neglecting all equations which involve modes
$\widehat{u}_{m,n}$ for large enough $|n|$, since there is an infinite
coupling (indeed, for each $n\in\Z$, the equation for $n$-th mode is
coupled with those of the two closest modulo 2 $n$-modes, thus
providing a tridiagonal system) and such a rough simplification would
lead to a system without solutions.  The main idea in the present
study is that higher frequencies are asymptotically small, since they
are terms of a convergent Fourier series. The fact that the series
converges derives from the abstract results in~\cite{Bei2005c,GR2005},
which can be even simplified in the case of an assigned smooth and
periodic pressure gradient. Hence, by assuming that $\widehat{f}_m$
vanishes fast enough as $|m|\to+\infty$ (and this is more than
reasonable for realistic flows, see also Section~\ref{sec_numtest}),
the assumed smallness of the higher modes is fully justified.
Therefore, once chosen a cut-off index $N\in\N$ ($N\ge2$ to keep at
least the first equation having zero as right-hand-side), we neglect
$\widehat{u}_{m,n}$ for all $n\in \Z$ such that $|n|>2N$ since, in
light of the aforementioned assumption, this does not affect the
solution in a significant manner.  From a practical viewpoint and
based on the symmetries highlighted above, this implies to drop off
$\widehat{u}_{m,2N+2}$ in the differential equation satisfied by
$\widehat{u}_{m,2N}$, and to only solve the equations for
$\widehat{u}_{m,n}$ with $n=0,2,\ldots,2N$.  In this way, for each
$m\in\N\cup\{0\}$, we end up with a closed tridiagonal system of
linear ordinary differential equations: it naturally represents the
differential counterpart of the tridiagonal matrix which is obtained
when evaluating the eigenvalues of the Mathieu functions by truncating
the corresponding algebraic linear system (see
e.g.~\cite{Alh1996,GPK2008}).

Before stating the tridiagonal system we actually address, we remark
that, by virtue of linearity, it is convenient for our purposes to
directly assume $\widehat{\sigma}_m=1$.  Moreover, we denote the
unknown by $\widehat{v}_{m,n}$ in order to stress that, besides
adopting a unitary pressure gradient, we are approximating
$\widehat{u}_{m,n}$ by truncation. In light of the above arguments and
derivations, for each $m\in\N\cup\{0\}$ we have to solve the following
boundary value problem:

\bigskip
{\small
  \begin{equation}
  \label{eq: ODE system}
  \left\{
    \begin{array}{rcl}
      \widehat{v}''_{m,0}
      -\bigg[\frac{ i \, \mathrm{Wo}^2_{a,m}}{2} \cosh(2\eta)\bigg]\,\widehat{v}_{m,0}
      +\frac{ i \, \mathrm{Wo}^2_{a,m}}{2}
      \widehat{v}_{m,2}
      &=&- \frac{a^2}{2 \nu}\cosh(2\eta)
      \\
      \widehat{v}''_{m,2}
      -\bigg[      2^2 +\frac{ i \, \mathrm{Wo}^2_{a,m}}{2}\cosh(2\eta)\bigg]\,\widehat{v}_{m,2}
      +\frac{ i \, \mathrm{Wo}^2_{a,m}}{4}(\widehat{v}_{m,4} +\widehat{v}_{m,0})
      &=&\frac{a^2}{4 \nu}
      \\
      \widehat{v}''_{m,4}
      -\bigg[      4^2 +\frac{ i \, \mathrm{Wo}^2_{a,m}}{2}\cosh(2\eta)\bigg]\,\widehat{v}_{m,4}
      +\frac{ i \, \mathrm{Wo}^2_{a,m}}{4}(\widehat{v}_{m,6} +\widehat{v}_{m,2})
      &=&\ \ 0
      \\
      \vdots\hspace{2cm}&=&\quad\vdots
      \\
      \widehat{v}''_{m,2n}
      -\bigg[      (2n)^2 +\frac{ i \, \mathrm{Wo}^2_{a,m}}{2}\cosh(2\eta)\bigg]\,\widehat{v}_{m,2n}
      +\frac{ i \, \mathrm{Wo}^2_{a,m}}{4}(\widehat{v}_{m,2n+2}+\widehat{v}_{m,2n-2})
      &=&\ \ 0
      \\
      \vdots\hspace{2cm}&=&\quad\vdots
      \\
      \widehat{v}''_{m,2N-2}
      -\bigg[      (2N-2)^2 +\frac{ i \, \mathrm{Wo}^2_{a,m}}{2}\cosh(2\eta)\bigg]\,\widehat{v}_{m,2N-2}
      +\frac{ i \, \mathrm{Wo}^2_{a,m}}{4}(\widehat{v}_{m,2N}+\widehat{v}_{m,2N-4})
      &=&\ \ 0
      \\
      \widehat{v}''_{m,2N}
      -\bigg[      (2N)^2 +\frac{ i \, \mathrm{Wo}^2_{a,m}}{2}\cosh(2\eta)\bigg]\,\widehat{v}_{m,2N}
      +\frac{ i \, \mathrm{Wo}^2_{a,m}}{4}(\widehat{v}_{m,2N-2})
      &=&\ \ 0
      \\
      \\
      \widehat{v}_{m,n}(b)=0\qquad  n=0,2,\dots,2N,&&
      \\
      \widehat{v}'_{m,n}(0)=0\qquad n=0,2,\dots,2N.&&
    \end{array}
\right.
\end{equation}
}
%
%
It is worth remarking how, instead of truncating an infinite
dimensional algebraic linear system in order to evaluate corresponding
eigenvalues (as in~\cite{GPK2008}), we performed a Galerkin-type
truncation of the system of ordinary differential equations which
stems from the Fourier separation of variables. It is in this sense
that we only \textit{indirectly} referred to Mathieu functions, thus
circumventing the numerical challenges associated therewith (indeed,
we end up with a very simple tridiagonal system, which can be
straightforwardly solved).

\begin{remark}
  An alternative pathway for solving the considered direct problem
  would be that of~\cite{GR2005}: By considering real (i.e. sine and
  cosine) Fourier series functions, to solve the following system for
  the unknowns $(\phi_m,\psi_m)$: $-n\psi_m=\nu\Delta\phi_m+1$,
  $n\phi_m=\nu\Delta\psi_m$, coupled with vanishing Dirichlet
  conditions.  Indeed, the above system can be used to prove the
  existence of time-periodic solutions in an abstract way. However, in
  view of our main objective of obtaining a benchmark solution for the
  flow in elliptical cross-sections, it seems less convenient to
  follow such a strategy, since it involves the bi-Laplacian operator,
  which seems not to be easily manageable by
  separation of variables in elliptical coordinates.
\end{remark}

In a practical framework, it is necessary to only solve the
system~\eqref{eq: ODE system} for $m\in\mathcal{M}^\star$, where
$\mathcal{M}^\star=\{0,1,\ldots,M^\star\}$, $M^\star$ denoting a sort
of upper bound for the harmonic content of the given flow rate
(i.e. its highest yet still relevant $t$-frequency is assumed to be
$\omega_{M^\star}$).  Moreover, once fixed the fluid (i.e. $\nu$) and
the size of the cross-section (i.e. $a$ and $b$), the solution
of~\eqref{eq: ODE system} is only affected by the adopted cut-off
index $N$.  Let us stress this point by introducing the more
descriptive notation $\widehat{v}^{(2N)}_{m,n}$, as an alternative to
$\widehat{v}_{m,n}$.  Hence, it is mandatory to choose $N$ large
enough to get a proper approximation for all the modes
$\widehat{v}^{(2N)}_{m,n}$, in particular for $m$ up to $M^\star$.
Let us denote by $N^\star$ the sought value of $N$, to remind that it
also depends on $M^\star$.  Two assets can be introduced to the
purpose: accuracy (magnitude of the cut modes must be negligible
compared to that of the kept ones) and independence from $N$
(truncation must only affect the kept modes in a negligible way).  As
discussed above, higher frequencies are expected to be small than
lower ones (by Fourier convergence); hence the accuracy asset can be
formulated by introducing the following non-dimensional quantity:
\begin{equation}
  \mu(m,N) = \frac{\|\widehat{v}_{m,2N}^{(2N)}\|_{\infty}}
  {\|\widehat{v}_{m,0}^{(2N)}\|_{\infty}},
\label{eq_mu_determ_Nstar}
\end{equation}
where $\|\widehat{v}\|_{\infty}=\max_{\eta\in E^\prime_\eta}
|\widehat{v}(\eta)|$, $E^\prime_\eta=[0,b]$ denotes the considered
$\eta$-range, and the $n=0$ mode is used to get a reference
value. (Such a reference was used for the numerical tests in
Section~\ref{sec_numtest}, turning out to be non-vanishing; otherwise,
alternative choices can be introduced by exploiting larger $n$
values.) In particular, once chosen a threshold $\bar{\mu}$, let us
define an index $N^\star_{\bar{\mu}}$ as the smallest integer $N$ such
that $\mu(m,N) \le \bar{\mu}$, for all $m\in\mathcal{M}^\star$. Hence,
by choosing $N\ge N^\star_{\bar{\mu}}$ we are guaranteed that
truncation only cuts negligible $n$-modes, for all relevant values of
$m$.  As regards independence from $N$, the following non-dimensional
quantity is introduced:
\begin{equation}
s(m,N)=\displaystyle \max_{n\in\mathcal{N}}
       \frac{\|\widehat{v}_{m,n}^{(2N)}-\widehat{v}_{m,n}^{(2N+2)}\|_{\infty}}
            {\|\widehat{v}_{m,0}^{(2N)}\|_{\infty}},
\label{eq_s_determ_Nstar}
\end{equation}
where $\mathcal{N}=\{0,2,\dots,2N-2\}$ (index $N$ is not considered,
for obvious reasons).  Then, once chosen a threshold $\bar{s}$, let us
define an index $N^\star_{\bar{s}}$ as the smallest integer $N$ such
that $s(m,N) \le \bar{s}$, for all $m\in\mathcal{M}^\star$.  Hence, by
choosing $N\ge N^\star_{\bar{s}}$ we are guaranteed that truncation
only negligibly affects the kept $n$-modes, for all relevant values of
$m$.  Finally, we naturally define
$N^\star=\max\left(N^\star_{\bar{\mu}},N^\star_{\bar{s}}\right)$, so
as to achieve both targeted assets.
%
%
\subsection{Flow in the ellipse: the inverse problem}
\label{sec_inv_ellisse_unsteady}
The interest in a benchmark solution for the inverse problem in
elliptical cross-sections has been highlighted in
Section~\ref{sec_intro}.  Nevertheless, some relevant results have
been only recently provided in~\cite{GPK2008}, grounded on the
solution of the direct problem by direct computation of the Mathieu
functions, as discussed in Section~\ref{sez_probdir_aux} (earlier
attempts in~\cite{RKHM2001} seem to be loosely grounded on a
theoretical basis).  In order to circumvent the aforementioned
challenges related to such a computation, we propose an alternative
method which is based on the fast-Fourier approach introduced in
Section~\ref{sez_probdir_aux}.  Thanks to linearity, we simply need to
link the Fourier coefficients of the flow rate with those of the
pressure gradient (as in the stationary case, thus
generalizing~\eqref{eq_f2lambda_staz}). The relevant solution strategy
is described below.

The flow rate associated with the
\textit{ansatz}~\eqref{eq_ansatz_doppiof} reads
\begin{equation*}
  f(t)=\sum_{m\in\Z} \widehat{f}_m \,e^{i\omega_m t} =
  \sum_{m\in\Z}
  \left[ \,
    \sum_{n\in\Z}
    \int_0^b \widehat{u}_{m,n}(\eta)
    \left( \int_0^{2\pi} e^{in\theta} \,J(\eta,\theta) \, d\theta \right) d\eta
  \right] \, e^{i\omega_m t}.
\end{equation*}
Then, once substituted $J$ from~\eqref{eq_J_etatheta}, by explicit
calculations we get
\begin{equation*}
  \frac{1}{a^2}
  \int_0^{2\pi} e^{in\theta} \, J(\eta,\theta) \, d\theta =
  \left\{
	  \begin{aligned}
	  & -\pi/2         & n=\pm2,
	  \\
	  & \pi\cosh(2\eta) & n=0,
	  \\
	  &  0              & n\in\Z\backslash\{0,\pm2\},
	  \end{aligned}
  \right.
\end{equation*}
so that the following relation can be easily obtained, for each $m\in\Z$:
\begin{equation}
  \widehat{f}_m = \frac{a^2\pi}{2}
  \int_0^b
  \Big(-\widehat{u}_{m,-2}(\eta)+2\cosh(2\eta)\,\widehat{u}_{m,0}(\eta)-\widehat{u}_{m,2}(\eta)\Big) 
  \,d\eta. 
  \label{fourier_f_u}
\end{equation}
Clearly, without an explicit knowledge of $\widehat{u}_{m,n}(\eta)$
(i.e. of the sought solution), it is not possible to evaluate the flux
by means of~\eqref{fourier_f_u}.  Nevertheless, by recalling from
Section~\ref{sez_probdir_aux} the approximate solution
$\widehat{v}_{m,n}$ (associated with a unitary pressure gradient), it
is straightforward to exploit~\eqref{fourier_f_u} in order to
approximate the Fourier coefficients $\widehat{\lambda}_m$ (of the
unknown pressure gradient) in terms of $\widehat{f}_m$.  Indeed,
thanks to linearity and by also recalling~\eqref{eq:symmetries}, the
following relation is immediately obtained, for each
$m\in\N\cup\{0\}$:
\begin{equation}
\label{eq:1}
  \widehat{\lambda}_m \cong \widehat{f}_m \, \frac{1}{\pi a^2}
  \left[\int_0^b\big(\cosh(2\eta)\,\widehat{v}_{m,0}(\eta)
                     -\widehat{v}_{m,2}(\eta)\big)\,d\eta\right]^{-1}.
\end{equation}
Latter expression deserves a remark: the bracketed quantity must be
non-vanishing for such an expression to be meaningful. If we consider
the exact solution $\widehat{u}_{m,n}$ of the infinite dimensional
system (based on explicit knowledge of the asymptotic behavior of
appropriate Hilbertian norms of eigenfunctions of the bi-Laplace
operator in the domain under consideration, as for the existence
results proved in~\cite{Bei2005c,GR2005}), then the map between
$\widehat{f}_m$ and $\widehat{\lambda}_m$ is one-to-one. This, in
turn, implies a non-vanishing denominator in the expression
corresponding to~\eqref{eq:1}.  Yet such a condition is not perfectly
guaranteed when considering an approximate solution
$\widehat{v}_{m,n}$.  Nonetheless, if the numerical approximation is
accurate enough (i.e. for a cut-off index $N$ large enough), the
approximate denominator is close enough to the exact one for the
non-vanishing condition to be fulfilled.  It must be pointed out that
the same potentially critical issue occurs when directly using the
Mathieu functions.  For instance, in~\cite{GPK2008} the authors
implicitly rely on a well-defined mapping between flow rate and
pressure coefficients for obtaining their solution. In that case, the
degree of the involved Mathieu functions must be large enough for the
associated numerical method to be well-posed, in complete analogy with
the present case.

Having introduced the main ingredients, it is now possible to
summarize the basic steps we propose for solving the inverse problem:
\begin{itemize}
\item[(S0)] We start by choosing a fluid (thus the kinematic viscosity
  $\nu$) and the size (namely $a$ and $b$) of the considered vessel
  cross-section.  We also assume a \textit{given} flow rate $f(t)$
  having period $T$;
\item[(S1)] We determine an integer $M$ such that the Fourier spectrum
  of $f(t)$ is suitably approximated by considering $M$ modes; popular
  metrics (e.g. the Pearson correlation coefficient used
  in~\cite{RKHM2001,GPK2008}) can be used to the purpose.  Practical
  constraints, due for instance to the experimental procedure
  implemented for measuring $f$, may come into play at this stage (for
  instance, the Nyquist limit is considered in~\cite{GPK2008}).
  Anyway, by the end of this step the Fourier components
  $\widehat{f}_m$, for $m=0,1,\ldots,M$, are known;
\item[(S2)] We fix an index $M^\star$ such that $M^\star \ge M$ and,
  once fixed desired (small enough) thresholds $\bar{\mu}$ and
  $\bar{s}$, we determine $N^\star$ as detailed at the end of
  Section~\ref{sez_probdir_aux}.  By the end of this step, the
  (accurate and truncation-independent) computation of the modes
  $\widehat{v}_{m,n}$ is achieved, for $m=0,1,\ldots,M^\star$
  and $n=0,2,\ldots,2N^\star$;
\item[(S3)] By considering the modes $\widehat{v}_{m,n}$ with $m=0,1,\ldots,M$
  (and $n=0,2,\ldots,2N^\star$), we compute $\widehat{\lambda}_m$
  by means of~\eqref{eq:1}, by also exploiting the relevant coefficient $\widehat{f}_m$
  determined at step (S1).
  By the end of this step, the coefficients $\widehat{\lambda}_m$ are available,
  for $m=0,1,\ldots,M$;
\item[(S4)] We finally compute the sought approximate solution
  by replacing~\eqref{eq_ansatz_doppiof} with the following summation:
  \begin{equation}
          u(t,\eta,\theta) \cong
          \sum_{m\in\bar{\mathcal{M}}} \widehat{\lambda}_m \, \varphi_m(\eta,\theta) \, e^{i\omega_m t},
          \quad \mathrm{with} \quad 
          \varphi_m(\eta,\theta) =
          \sum_{n\in\bar{\mathcal{N}}^\star} \widehat{v}_{m,n}(\eta) \, e^{in\theta},
  \label{def_somma_phi}
	\end{equation}
  where $\bar{\mathcal{M}}=\{-M,-M+1,\ldots,M-1,M\}$ and
  $\bar{\mathcal{N}}^\star=\{-2N^\star,-2N^\star+2,\ldots,2N^\star-2,2N^\star\}$,
  having previously obtained the modes associated with negative values
  of $m$ and $n$ by conjugacy, as detailed in
  Section~\ref{sez_probdir_aux}.
\end{itemize}

\noindent
We conclude by adding two remarks on Step (S2).  Firstly, to
accomplish such a task, we need to solve the system~\eqref{eq: ODE
  system} several times, namely for several $N$. However, this burden
is not specifically introduced by our strategy: for instance, to
iterate computations up to reaching an adequate level of accuracy is
also compulsory when directly using the Mathieu functions (see
e.g.~\cite{GPK2008}).  Indeed, our strategy (which simply involves the
solution of tridiagonal systems) holds potential for a more
computationally-efficient procedure than direct management of the
Mathieu functions.  Latter remark regards the introduction of
$M^\star$.  In particular, if one is interested in a single flow rate
it makes sense to directly choose $M^\star=M$, so as to minimize
computations.  However, it can be of interest to also address multiple
flow rate conditions, e.g. to identify the effect of specific
$m$-harmonics on the solution (this is of interest, for instance, for
the development of pulsatility-based devices, as mentioned in
Section~\ref{sec_intro}).  In such a case, it is convenient to choose
$M^\star > M$ when taking step (S2), in particular by adopting an
$M^\star$ large enough to account for all the foreseen $m$-harmonic
contributions.  Indeed, once built the corresponding
$\widehat{v}_{m,n}$ ``basis'', it is possible to immediately explore
many solutions, by simply iterating steps (S3)-(S4) in correspondence
of any given set of $\widehat{f}_m$ coefficients.  This leads to a
very computationally-cheap procedure, since it is only required to
evaluate some integrals, besides trivial summations.  This aspect
further supports the effective exploitability of the solution strategy
we propose.

\section{Numerical results}
\label{sec_numtest}
The numerical strategy detailed in
Section~\ref{sec_inv_ellisse_unsteady} is hereafter exemplified, based
on measured data of blood and cerebrospinal fluid (CSF) flow under
physiological conditions.  As anticipated in Section~\ref{sec_intro},
despite the exploitation of physiological flow conditions, the
considered test-cases do not lay strong claims of being
physiologically representative, due to the assumption of a
fully-developed flow, which can be hardly fulfilled in real
situations.  Main aim of such simulations is to quantitatively assess
the gain in computational efficiency which is enabled by the proposed
Fourier-based numerical solution, as compared to more elaborate finite
element solvers. This, in turn, renders our benchmark solution
effectively usable for developing more ambitious numerical approaches,
up to serving as a debugging tool for complex 3D codes.

We point out that, besides the simulations described below, which are
associated with mild eccentricity values (above $0.5$), we
successfully manage to also apply our numerical strategy to more
idealized geometries with eccentricity e.g. below $0.1$, thus
succeeding where previous approaches encountered
difficulties~\cite{HZ1998}. However, corresponding numerical results
are very similar to those reported in the sequel and therefore they
are not reported, for ease of presentation.

\subsection{Blood flow in the internal carotid artery}
\label{sec_res_ICA}

\begin{figure}[b!]
  \begin{center}
  {\includegraphics[width=15.9cm]{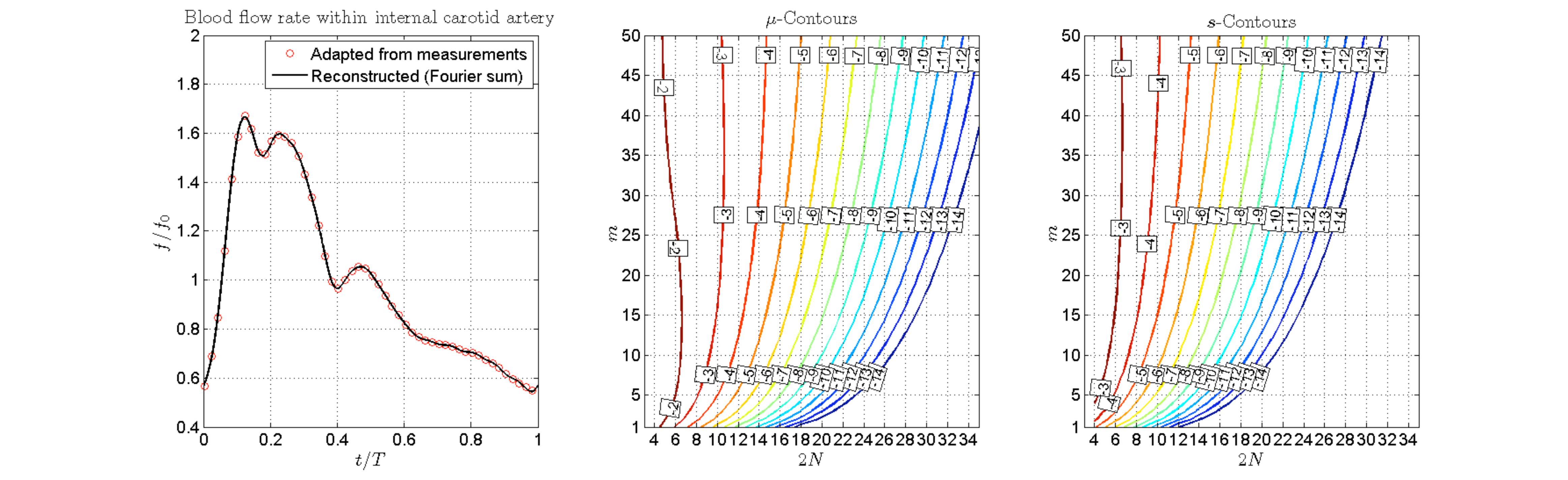}}
  \caption{Blood flow within ICA.  (Left) Normalized flow rate
    $f/f_0$: data, adapted from \cite{Hoi2010}, are plotted against
    Fourier reconstruction based on $M=15$ modes.  (Middle) Contour
    plots for $\mu(m,N)$: e.g. it is necessary to choose $2N=34$ for
    the relative magnitude of discarded modes to be below $10^{-12}$,
    for all the considered $m$ values. (Right) Contour plots for
    $s(m,N)$: e.g. it is necessary to choose $2N=28$ to get a relative
    sensitivity to truncation below $10^{-12}$, for all the considered
    $m$ values.}
\label{fig_flusso_ICA_Nstar}
\end{center}
\end{figure}

\begin{figure}[h!]
\begin{center}
{\includegraphics[width=15.9cm]{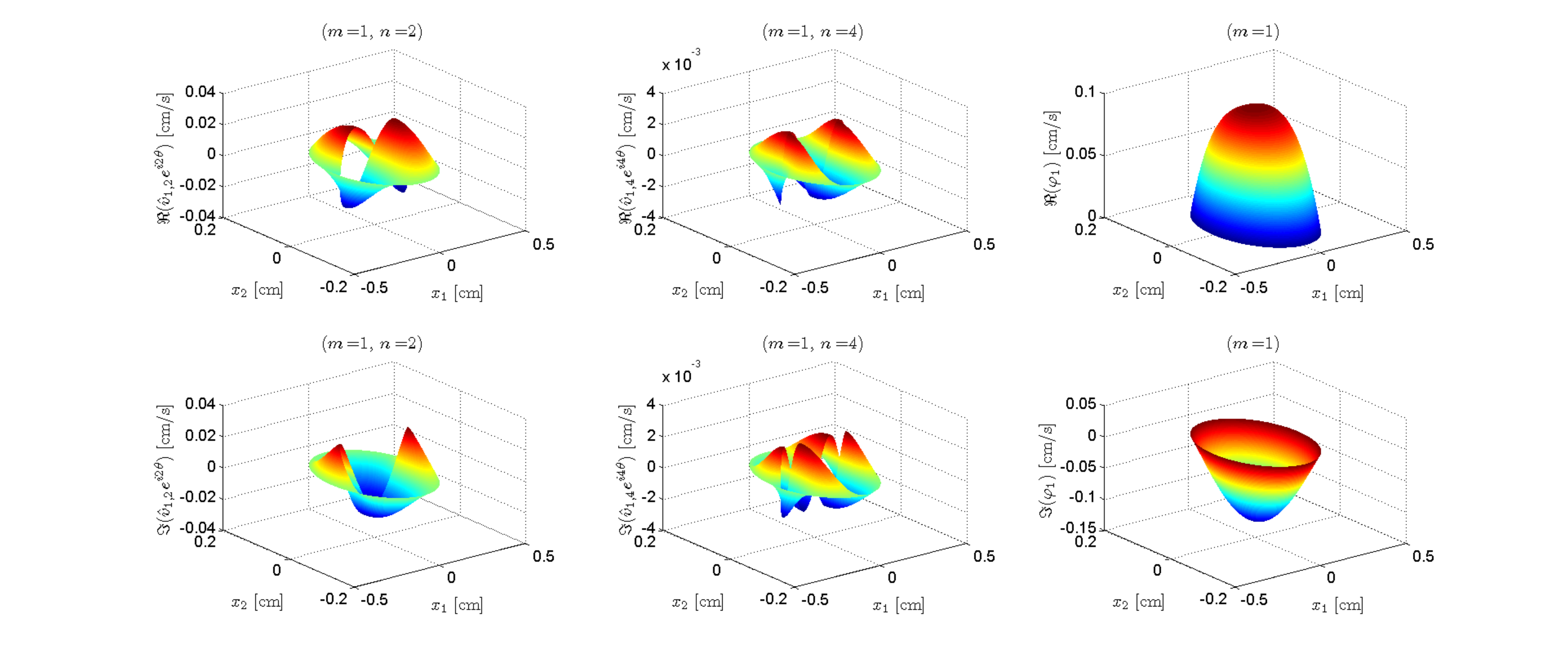}}
\caption{Blood flow within ICA.  Real part (upper row) and imaginary
  part (lower row) of selected solution components, computed with
  $2N=34$: $\widehat{v}_{1,2}(\eta)e^{i 2\theta}$ (left column),
  $\widehat{v}_{1,4}(\eta)e^{i 4\theta}$ (middle column) and
  $\varphi_1(\eta,\theta)$ (right column). Latter quantity, defined
  in~\eqref{def_somma_phi}, accounts in particular for all the $m=1$
  modes.}
\label{fig_ICA_componenti_m1}
\end{center}
\end{figure}

A flow rate waveform for blood flow within the human internal carotid
artery (ICA) was adapted from~\cite{Hoi2010}, see
Fig~\ref{fig_flusso_ICA_Nstar} (left); corresponding period and
period-averaged flow rate are $T=0.95$~[s] and $f_0=4.11$~[cm$^3$/s],
respectively.
An average radius for such a vessel is
$0.25$~cm~\cite{Krejza2006}. Hence, by assuming e.g. an eccentricity
equal to $0.6$, we consistently introduced an elliptical cross-section
with semi-axes $\alpha=0.25$~[cm] and $\beta=0.15$~[cm], being the
corresponding $\eta$-range $E'_{\eta}=[0,b]$, with $b=0.69$ (a more
realistic vessel deformation, e.g. due to compression by the
surrounding tissue, might imply some larger value for $\alpha$;
however such a refined estimate is beyond the scope of the present
simulations).  Furthermore, a characteristic diameter for the
considered section is $\tilde{\delta}=\alpha+\beta=0.4$~[cm], while
the section-averaged speeds is $\bar{w}=f_0/A=34.9$~[cm/s], where
$A=\pi \alpha \beta$ represents the cross-sectional area.  The adopted
flux was then approximated by $M=15$ modes; corresponding Fourier sum,
also reported in Fig.~\ref{fig_flusso_ICA_Nstar}, suitably replaces
the data at hand (associated Pearson correlation coefficient differs
from $1$ by less than $10^{-3}$). Moreover, once introduced a
characteristic speed $\tilde{w}=A^{-1}
\max_{t/T\in[0,1]}f(t)=58.36$~[cm/s] and by assuming $\nu=3.5 \cdot
10^{-2}$~[cm$^2$/s]~\cite{rogers2010blood}, it is possible to label
the considered blood flow with a Reynolds number
$Re=\tilde{w}\tilde{\delta}/\nu \cong 670$.  In addition, once defined
a characteristic frequency $\widetilde{\omega}=\left(\sum_{m=0}^{M}
  |\widehat{f}_m| \, \omega_m\right)/ \left(\sum_{m=0}^{M}
  |\widehat{f}_m|\right)$, it is possible to also introduce a
characteristic Womersley number, namely
$\mathrm{Wo}=(\tilde{\delta}/2)\sqrt{\widetilde{\omega}/\nu} \cong
1.5$.  Despite the scarcity of results on stability for pulsatile
flows (which seem to be often contradictory with one another already
for the circular cross-section~\cite{Trip2012}), the above figures
suggest that the flow at hand is laminar.  This was derived by
considering available experimental thresholds for the circular
cross-section~\cite{Carpinlioglu2001,Haddad2010}, possibly defined
through empirical fittings as in~\cite{Peacock1998} (to the purpose,
it may be appropriate to also define relevant Reynolds and Womersley
numbers in a slightly different way~\cite{Carpinlioglu2001}, yet we
verified that this does not affect the above statement on stability
for the considered flow).  That said, more detailed considerations on
laminar stability cannot be drawn in the present scope: they certainly
require more consolidated acquisitions, as expressly remarked
e.g. in~\cite{Carpinlioglu2001}.  Step (S2) of the procedure proposed
in Section~\ref{sec_inv_ellisse_unsteady} was then addressed, by
choosing $M^\star=50$ for the sake of illustration.  More in detail,
the contours shown in Fig.~\ref{fig_flusso_ICA_Nstar} were firstly
derived for both $\mu(m,N)$ and $s(m,N)$, respectively defined
in~\eqref{eq_mu_determ_Nstar} and~\eqref{eq_s_determ_Nstar}, for
$m=1,\ldots,M^\star$ and $N$ large enough. To the purpose, we coded
the boundary value problem~\eqref{eq: ODE system} within Matlab
environment and we solved it by a shooting technique.
As a result, by choosing e.g. $\bar{\mu}=\bar{s}=10^{-12}$, it is
immediate to get $N^\star_{\bar{\mu}}=17$ and $N^\star_{\bar{s}}=14$
from the relevant plots in Fig.~\ref{fig_flusso_ICA_Nstar}, so as to
finally determine $N^\star=17$.  Please also notice how, for any fixed
$m$, $\| \widehat{v}_{m,n} \|_{\infty} \geq \|
\widehat{v}_{m,n+2}\|_{\infty}$ for $n=0,...,2N-2$, thus confirming
the expected relative less influence of higher $n$-modes.  Relevant
modes $\widehat{v}_{m,n}$ were thus obtained by solving~\eqref{eq: ODE
  system} with $N=N^\star$; Fig~\ref{fig_ICA_componenti_m1} shows some
solution components, including $\varphi_1$ defined
in~\eqref{def_somma_phi}.

\bigskip
\bigskip

\begin{figure}[t!]
\begin{center}
  {\includegraphics[width=15.9cm]{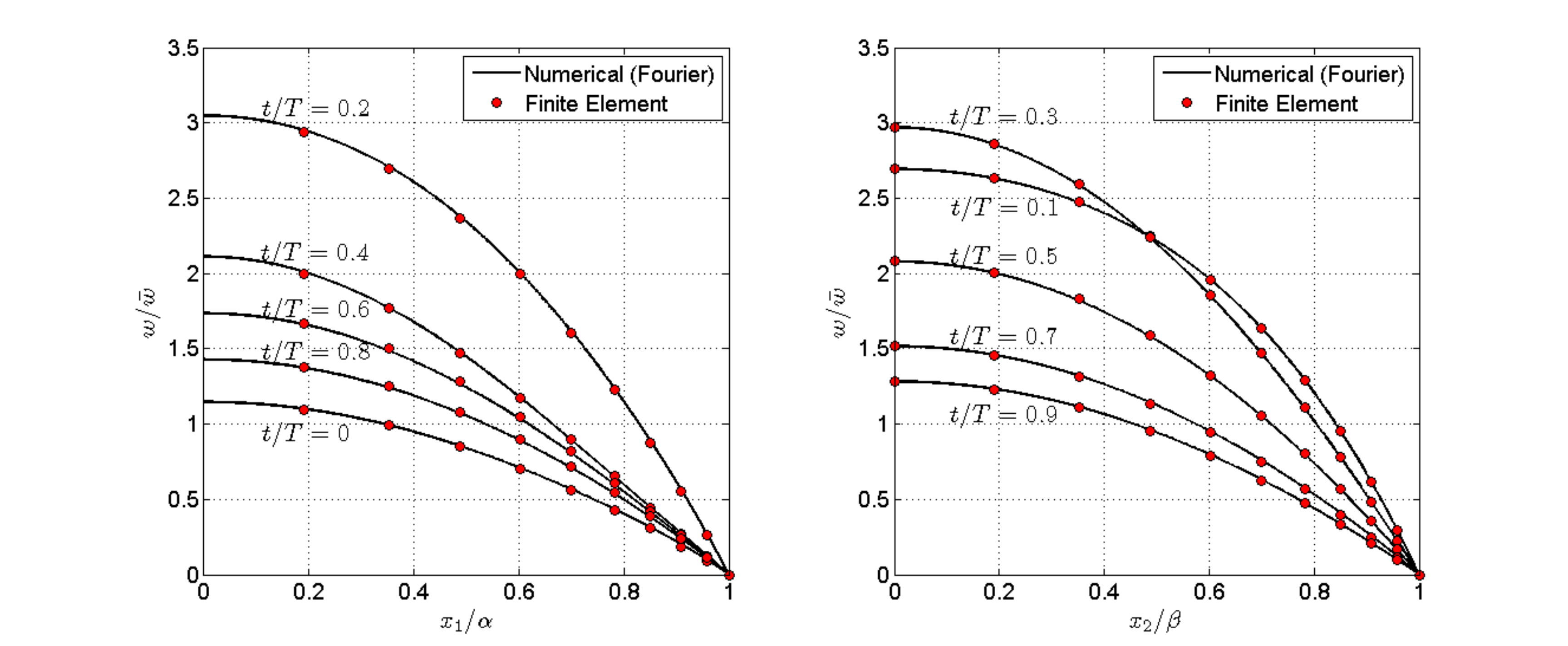}}
\caption{Blood flow within ICA. Normalized velocity profiles $w/\bar{w}$
along the major (left) and minor (right) semi-axes of the ellipse, for selected
non-dimensional times $t/T$. Results obtained by the proposed Fourier-based
approach (solid curves) are compared with those achieved by a commercial
FE solver (filled circles).}
\label{fig_ICA_matlab_adina}
\end{center}
\end{figure}

We then compared our numerical results with those achieved by means of
the commercial finite element (FE) solver ADINA 8.8.1 (ADINA R$\&$D
Inc., MA, USA), available to the group.  Such a solver uses a standard
Galerkin formulation (stabilized by upwinding for higher Reynolds
numbers~\cite{Bathe1995193}), while time-advancing is implemented by a
first-order Euler backward method.  As for FE simulations, a pipe with
length $\ell = 160 \tilde{\delta}$ was defined, in order to obtain a
fully-developed flow in the central portion of the domain (such a
condition was \textit{a posteriori} checked).  Furthermore, the flow
rate shown in Fig~\ref{fig_flusso_ICA_Nstar} was imposed (by software
coding) at the inlet cross-section, the no-slip (i.e. Dirichlet)
boundary condition was enforced on vessel wall, and a reference
pressure value was imposed at the outlet cross-section (such value is
immaterial, due to the incompressible formulation).  In addition, as
regards simulation set-up, both space- and time-discretization were
incrementally refined, up to obtain discretization-independent
results.  In particular, pipe domain was discretized by nearly
$8.2\cdot10^{5}$ second-order accurate brick elements, and $7$
pulsation periods were simulated (time-periodicity was obtained after
$4$ periods).  Finally, all simulations were run on a single core of a
PC with Intel Core i7-960 3.20 GHz CPU and 24 GB RAM.  Obtained
results are shown in Fig~\ref{fig_ICA_matlab_adina}; as expected, a
very good agreement is achieved between the considered approaches,
main difference residing in computational times.  Indeed, once
accomplished in roughly $3$ minutes the preliminary step (S2) (see
Section~\ref{sec_inv_ellisse_unsteady}) so as to get
truncation-independent results, the computation is finalized in a very
few seconds.  Conversely, runtime for the FE simulation was about $9$
days; additional days were necessary to set-up the simulation (namely
to get grid-independence, periodicity and fully-developed flow
conditions), so that the overall computational time for the FE run can
be estimated over $20$ days. Corresponding speed-up factor is close to
$10^{4}$.

\subsection{Cerebrospinal fluid flow in the spine}
\label{sec_res_CSF}

\begin{figure}[b!]
  \begin{center}
  {\includegraphics[width=15.9cm]{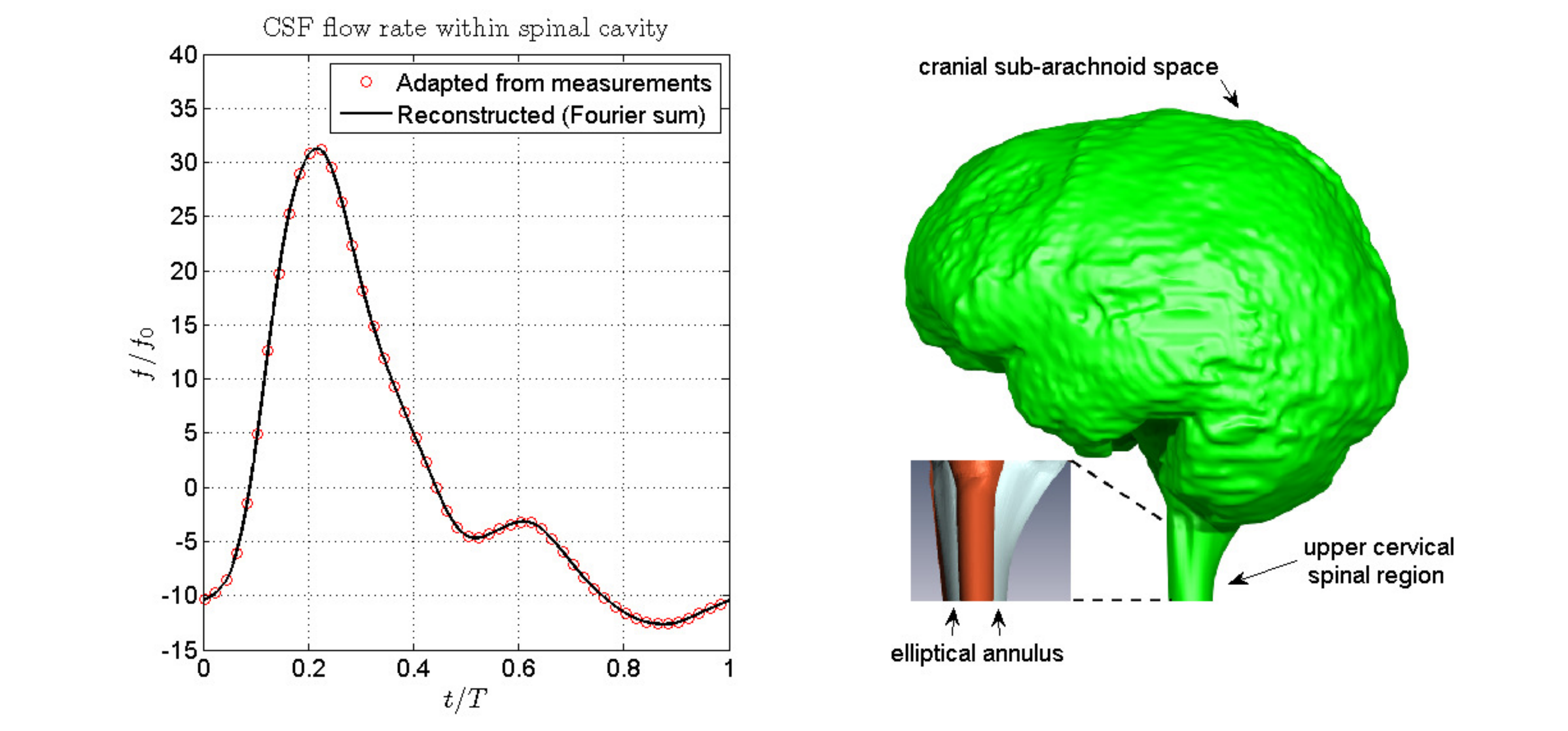}}
  \caption{CSF flow in the spine. 
  (Left) Normalized flow rate $f/f_0$: data, adapted from \cite{GPK2008}, are plotted against
  Fourier reconstruction based on $M=15$ modes.
  (Right) 3D Reconstruction of the CSF domain, as obtained from
  patient-specific magnetic resonance imaging (MRI). 
  Cranial sub-arachnoid space~\cite{Irani2008} and upper cervical spinal region are shown;
  location of the considered annulus is shown in the inset.}
\label{fig_flusso_CSF_3Drec}
\end{center}
\end{figure}
\begin{figure}[t!]
\begin{center}
{\includegraphics[width=15.9cm]{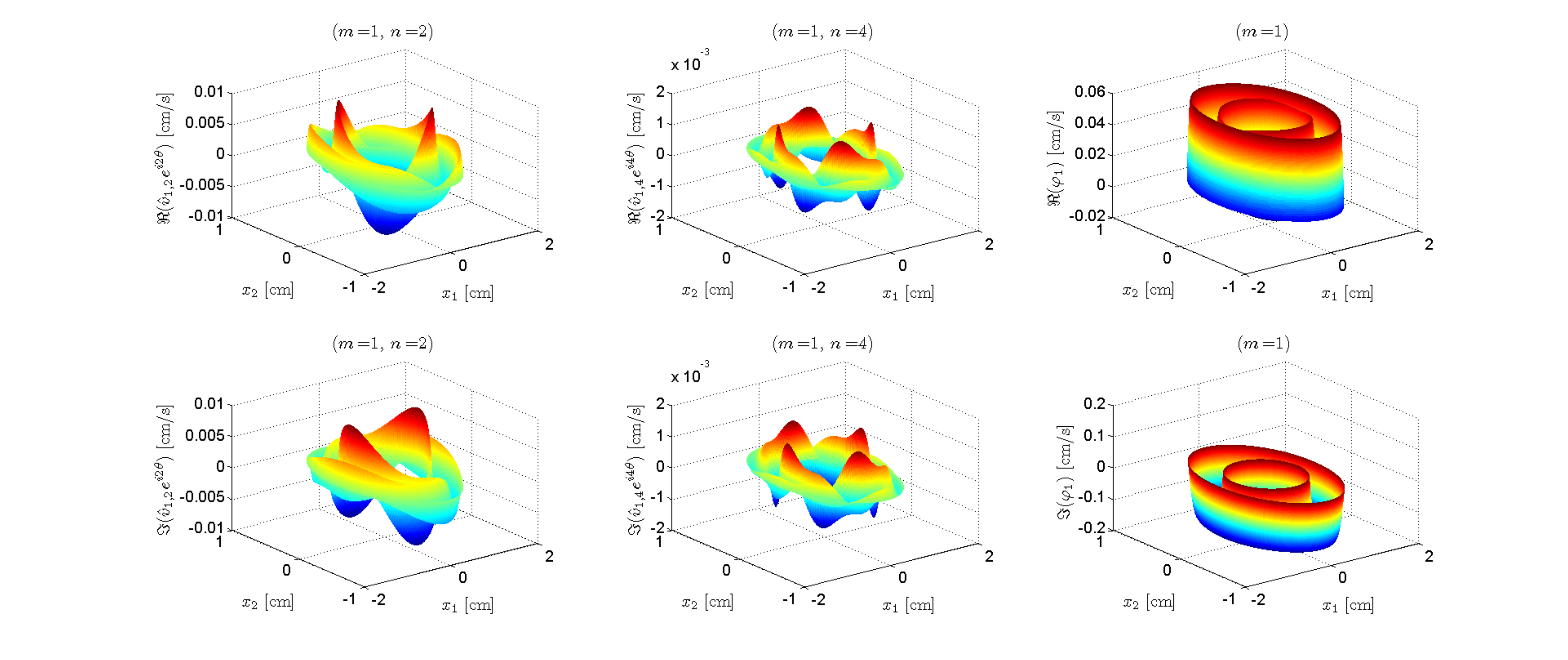}}
\caption{CSF flow in the spine.
 Real part (upper row) and imaginary part (lower row) of selected solution components,
 computed with $2N=34$:
 $\widehat{v}_{1,2}(\eta)e^{i 2\theta}$ (left column),
 $\widehat{v}_{1,4}(\eta)e^{i 4\theta}$ (middle column)
 and $\varphi_1(\eta,\theta)$ (right column). 
 Latter quantity, defined in~\eqref{def_somma_phi},
 accounts in particular for all the $m=1$ modes.}
\label{fig_CNS_componenti_m1}
\end{center}
\end{figure}

A flow rate waveform for cerebrospinal fluid (CSF) in the upper cervical region of
the human spine~\cite{Irani2008}) was adapted
from~\cite{GPK2008}, see Fig~\ref{fig_flusso_CSF_3Drec};
corresponding period-averaged flow rate is $f_0 = -0.11$~[cm$^3$/s] (negative sign here
indicates that it is directed towards the lumbar region).
Corresponding period was assumed equal to that one of the cardiac cycle, reported in
Section~\ref{sec_res_ICA}.
The considered cross-section can be approximated by the
annulus between two confocal ellipses with $\alpha_2=1.11$~[cm],
$\beta_2=0.93$~[cm] (from which $a=0.61$~[cm]) and
$\beta_1 = 0.43$~[cm],
being the corresponding $\eta$-range $E'_{\eta}=[b_1,b_2]$, with
$b_1=0.66$ and $b_2=1.21$.  Furthermore, a characteristic thickness
for the considered section is
$\tilde{\tau}=\left(\alpha_2-\alpha_1+\beta_2-\beta_1\right)/2=0.43$~[cm],
while the section-averaged speeds is $\bar{w}=f_0/A=-0.47$~[cm/s],
where $A$ represents the cross-sectional area.  Such a flux was then
approximated by $M=15$ modes; corresponding Fourier sum, also reported
in Fig.~\ref{fig_flusso_CSF_3Drec}, suitably replaces the data at hand
(associated Pearson correlation coefficient differs from $1$ by less
than $10^{-6}$).  Moreover, by reasoning as in
Section~\ref{sec_res_ICA} with $\tilde{\tau}$ in place of
$\tilde{\delta}$, and by assuming
$\nu=10^{-2}$~[cm$^2$/s]~\cite{LYA2001}, it is possible to label the
considered blood flow with a Reynolds number
$Re=\tilde{w}\tilde{\tau}/\nu \cong 64$ and a Womersley number
$\mathrm{Wo}=(\tilde{\tau}/2)\sqrt{\widetilde{\omega}/\nu} \cong 5.5$.
In light of the relevant references introduced in
Section~\ref{sec_res_ICA}, the above figures suggest that the flow at
hand is laminar.  However, also in this case more consolidated results
are needed for a stronger statement on stability (available data for
the circular cross-section~\cite{Trip2012,Carpinlioglu2001} better
describe almost purely oscillating flows like the present one, yet not
in an annular section); such a study is clearly beyond the present
scope.  By proceeding as for the ICA, we then determined $N^\star$ by
considering $M^\star=50$ $m$-modes: once chosen
$\bar{\mu}=\bar{s}=10^{-12}$ we got $N^\star=17$ (associated contours
are not shown for brevity, being similar to those in
Fig~\ref{fig_flusso_ICA_Nstar}).  Relevant modes $\widehat{v}_{m,n}$
were thus obtained by solving~\eqref{eq: ODE system} with $N=N^\star$;
Fig~\ref{fig_CNS_componenti_m1} shows some solution components,
including $\varphi_1$ defined in~\eqref{def_somma_phi}.

\begin{figure}[t!]
\begin{center}
  {\includegraphics[width=15.9cm]{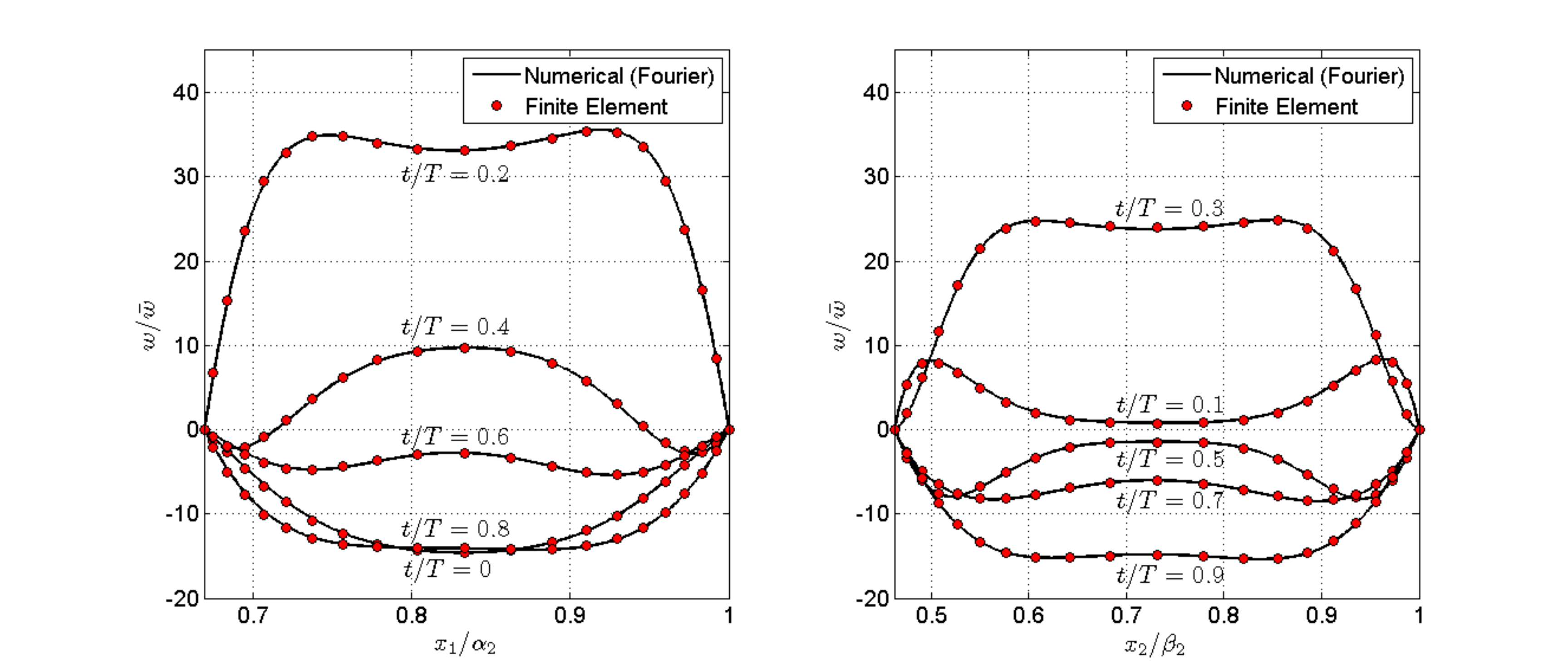}}
\caption{CSF flow in the spine. Normalized velocity profiles $w/\bar{w}$
along the major (left) and minor (right) semi-axes of the annulus, for selected
non-dimensional times $t/T$. Results obtained by the proposed Fourier-based
approach (solid curves) are compared with those achieved by a commercial
FE solver (filled circles).}
\label{fig_CSF_matlab_adina}
\end{center}
\end{figure}

As for the ICA test-case, we compared our numerical results with those
achieved by the commercial solver ADINA (relevant details can be
recalled from Section~\ref{sec_res_ICA}).  For such FE simulations, a
pipe with length $\ell \cong 100 \bar{\tau}$ was defined, in order to
obtain a fully-developed flow in the central portion of the domain
(such a condition was \textit{a posteriori} checked).  Furthermore,
the flow rate shown in Fig~\ref{fig_flusso_CSF_3Drec} was imposed (by
software coding) at the inlet cross-section, the no-slip
(i.e. Dirichlet) boundary condition was enforced on vessel walls, and
a reference pressure value was imposed at the outlet cross-section
(such value is immaterial, due to the incompressible formulation).  In
addition, as regards simulation set-up, both space- and
time-discretization were incrementally refined, up to obtain
discretization-independent results.  In particular, pipe domain was
discretized by nearly $3.8\cdot10^{5}$ second-order accurate brick
elements, and $7$ pulsation periods were simulated (time-periodicity
was obtained after $2$ periods).  Obtained results are shown in
Fig~\ref{fig_CSF_matlab_adina} (they are obviously identical to those
reported in~\cite{GPK2008}). As expected,
a very good agreement is achieved between the proposed Fourier-based
solution and that provided by the finite element solver, main
difference residing in computational times for the present case as
well.  Indeed, once accomplished in roughly $3$ minutes the
preliminary step (S2) (see Section~\ref{sec_inv_ellisse_unsteady}) so
as to get truncation-independent results, the computation is finalized
in a very few seconds.  Conversely, runtime for the FE simulation was
about $4$ days (less than for ICA thanks to the lower average speeds);
additional days were necessary to set-up the simulation (namely to get
grid-independence, periodicity and fully-developed flow conditions),
so that the overall computational time for the FE run can be estimated
over $10$ days. Corresponding speed-up factor is well over $10^{3}$
also in the present case.

\section{Concluding remarks}
\label{sec_conclusions}

We proposed a novel numerical method for solving the inverse problem
of pulsatile viscous flows in elliptical vessels and annuli.  In
particular, we addressed a fully-developed flow of an incompressible
Newtonian fluid, so as preliminarily derive some analytical relations
(otherwise hardly achievable) upon which to build our numerical
solution.  We are therefore aware of the fact that the solution we
propose is hardly applicable to 3D finite-length complex domains (such
a limitation, despite being obvious, is often understated by the many
authors who also assumed a fully-developed flow).  On regard, under
quite general assumptions (yet still adopting an incompressible
Newtonian fluid), many valuable approaches for problems with flow rate
conditions have been proposed in the literature (see
e.g.~\cite{FVV2010}), and they are the subject of ongoing research.
Moreover, a fully numerical approach is the only viable when also
considering non-Newtonian fluids, already when tackling stationary
problems (possibly within curved or non-straight vessels as
in~\cite{APS2007}).  Nevertheless, and in spite of the linear context
within which it was derived, our solution to the considered inverse
problem is original and non-trivial.  Indeed, it provides an easily
computable benchmark (as well as an approximation) for such a problem.
In particular, by invoking separation of variables (in elliptical
coordinates) we obtained a tridiagonal system of linear, ordinary
differential equations and we obtained the sought solution by means of
a Galerkin-type truncation.  This method seem to theoretically
corroborate the early attempts in~\cite{RKHM2001} (in which
proper references to the inverse problem are missing).  Moreover, it
represents an alternative to the numerical approach in~\cite{GPK2008},
heavily based on a wise computation strategy of the Mathieu functions.
Indeed, truncation of the aforementioned system represents the
differential counterpart of the truncation of the infinite dimensional
algebraic linear system in~\cite{GPK2008}, which is performed in order
to evaluate corresponding eigenvalues.  In this sense, we only
indirectly referred to Mathieu functions, thus circumventing the
numerical challenges associated therewith
(see~\cite{GPK2008,HZ1998,SW2009}).  As a consequence, even if the
considered differential system is not completely integrable (due to
the coupling between Fourier modes which is implied by the elliptical
cross-section already when addressing the stationary problem), our
solution can be easily computed, even for small ellipticity values.
Such an asset was shown through the numerical simulations reported in
Section~\ref{sec_numtest}, based on measured flow rates for blood in
the internal carotid artery (ICA), and cerebrospinal fluid (CSF) flow
in the upper cervical region of the human spine.  Indeed, despite
being expected, the main advantage of the proposed method, as compared
to fully numerical approaches, resides in computational efficiency: a
speed-up over $10^3$ was obtained for both the aforementioned
test-cases.  Such a gain is comparable with the one reported
in~\cite{GPK2008}, even if simulation set-up times seem to be slightly
understated therein (in particular, computation of the Mathieu
functions was recognized as the most critical part, and several
related convergence results are reported in such a valuable paper, yet
at the same time corresponding computational burden was judged as
inexpensive~\cite{GPK2008}).  An additional, potential advantage of
the proposed Fourier-based approach resides in the possibility of an
extensive and detailed control over numerical implementation (as
opposite to black-box libraries usually exploited for computing
special functions), which could be further optimized e.g. by using
consolidated techniques for Fast Fourier solvers~\cite{GL2004}.

These results effectively support the exploitation of the proposed
method in a wider research scope, e.g. as a debugging tool and/or an
improved source of boundary data for more ambitious numerical
approaches, possibly based on more realistic data. On regard and with
reference to the study of blood flow, the proposed solution could
provide an enhanced inflow condition for those approaches which merely
scale a (steady) parabolic solution by a factor accounting for flow
rate variability (see e.g.~\cite{GBBCFMS2011}).  Moreover, the
proposed solution can be directly used for computing wall shear stress
distribution along the adopted cross-section boundary.  Indeed, vessel
ellipticity affects such a stress~\cite{HZ1998}, which in turn
has been identified as a key player in the onset and development of relevant
pathologies likes atherosclerosis, even if
several acquisitions needs to be consolidated
(see e.g.~\cite{RKHM2001,Caro2009}). 
Furthermore, still within the context of
bio-fluid dynamics, many investigations are being carried out for
elucidating the role of CSF dynamics within human spine, in relation
to the onset and progression of relevant diseases like
syringomyelia.  The proposed benchmark flow holds potential for
application to such studies as well, in particular for assigning more
accurate inflow/initial conditions in the cervical portion of the
spinal fluid domain (which is usually described as a circular annulus between
two finite-length, compliant cylinders~\cite{B2009,CK2012,E2012};
wave propagation is usually studied therein).
It is worth mentioning that relevant investigations for such
fluid-structure interaction problems are
being performed in the context of hemodynamics as well (see e.g.~\cite{Canic2006}).
Furthermore, besides studies purely focusing on fluid dynamics, the
proposed solution can be effectively applied to additional biomedical
applications (as anticipated, the biomedical field is an elective
application domain, since flow rate is a main available physical
datum), e.g. to targeted delivery problems.  For instance, the effect
of flow pulsatility on magnetic particle targeting was originally
exemplified in~\cite{BS2011}, in order to highlight peculiar effects
with respect to the stationary case~\cite{Haverkort2009} widely
adopted when addressing particle targeting in blood capillaries. Such
an enhancement straightforwardly applies to magnetoresponsive carriers
at large, e.g. to the controlled navigation of MRI-guided microrobots
in the vasculature~\cite{Ferreira2011}. Indeed, there is a growing
interest in the robotic field for the development of
micro/miniaturized interventional tools to be navigated (e.g. through
magnetic fields) within bodily fluids. These systems, currently
defined as microrobots ~\cite{Nelson2010} despite their intrinsic
simplicity (system compartmentalization is indeed almost inapplicable
at the micro scale), hold potential for effective application to
medical fields, and both blood flow and CSF flow were identified as
elective carrying fluids. To further support this point, it is worth
mentioning that, besides carrier transport along the mainstream,
active microrobot navigation was also proposed, e.g. in the CSF
flow. More precisely, magnetic navigation of helical microrobots was
proposed for application in the spinal CSF~\cite{Abbott2010} (being
the helical shape inspired by the propulsion strategy of bacterial
flagella, well suited for low-Reynolds-number regimes).
Interestingly, relevant CSF velocity profile,
together with corresponding fluidic actions on the microrobot, seem to
be slightly overlooked in these studies  (while attention is payed to other physical
aspects, up to gravity, in spite of its minor effects at the micro
scale~\cite{Abbott2011}). This seems to be partly due to
the lack of easily computable - yet physically derived - flow
conditions; such an issue might be mitigated by our benchmark flow
solution, thus contributing to better assess the feasibility of such
long-term approaches in bodily fluids, even if still in a simplified
context. Finally, an additional, remarkable application field has been
recently identified, namely the one of energy harvesting from bodily
fluids. Indeed, use of magnetic fields was proposed as a key powering
strategy for most of the aforementioned microrobots, mainly in
response to the impossibility of using on-board batteries at the micro
scale (alternative approaches based on harnessing the power of
biological entities like bacteria were also proposed, yet their
discussion is out of the present scope).  In fact, powering (closely
linked to actuation) appears as the real bottleneck currently
hampering microrobotic applications, and energy harvesting holds
promise to solve such a major issue. In this spirit yet in the wider
context of implantable miniaturized electronic devices for biomedical
applications, an implantable fuel cell has been recently
proposed~\cite{Sarpeshkar2012}, powered through energy scavenging from
the CSF. Relevant estimates are reported in the cited reference in
order to assess powering effective sustainability (they address oxygen
and glucose transport within the CSF domain, mainly in the brain
region but also encompassing the spinal one). However, it is expressly
recognized in~\cite{Sarpeshkar2012} that a deeper insight into the
actual CSF dynamics is required to accurately assess the performance
of the proposed system.  Such a quest further supports the development
of modeling strategies for the CSF flow (yet such a point can be
immediately extended to other fluid conditions), and the benchmark
flow solution we propose, despite the inherent simplifications we
repeatedly pointed out above, could contribute to take a leap ahead in
many of the aforementioned applications areas.

To conclude, while
developing the proposed numerical method we deliberately sacrificed
some degree of physiological representativeness, yet in the name of
direct and efficient computability, and thus wider usability of the
obtained benchmark flow solution.  The many applications discussed
above encourage to exploit our solution in an wide interdisciplinary
context, so as to further push the corresponding scientific and
technological frontiers.

\section*{Acknowledgments}

The authors would like to thank Costanza Diversi and Byung-Jeon Kang for the
3D reconstruction reported in Fig~\ref{fig_flusso_CSF_3Drec}.


\begin{thebibliography}{10}

\bibitem{Alh1996}
{\sc F.~A. Alhargan}, {\em A complete method for the computations of {M}athieu
  characteristic numbers of integer orders}, SIAM Rev., 38 (1996),
  pp.~239--255.

\bibitem{APS2007}
{\sc N.~Arada, M.~Pires, and A.~Sequeira}, {\em Viscosity effects on flows of
  generalized {N}ewtonian fluids through curved pipes}, Comput. Math. Appl., 53
  (2007), pp.~625--646.

\bibitem{Ferreira2011}
{\sc L.~Arcese, M.~Fruchard, F.~Beyeler, A.~Ferreira, and B.~J. Nelson}, {\em
  Adaptive backstepping and {MEM}s force sensor for an {MRI}-guided microrobot
  in the vasculature}, in Proceedings of 2011 IEEE International Conference on
  Robotics and Automation, Shanghai, China, May 2011, pp.~4121--4126.

\bibitem{Bathe1995193}
{\sc K.~J. Bathe, H.~Zhang, and M.~H. Wang}, {\em Finite element analysis of
  incompressible and compressible fluid flows with free surfaces and structural
  interactions}, Comput. Struct., 56 (1995), pp.~193--213.

\bibitem{Bei2005c}
{\sc H.~{Beir{\~a}o da Veiga}}, {\em Time periodic solutions of the
  {N}avier-{S}tokes equations in unbounded cylindrical domains---{L}eray's
  problem for periodic flows}, Arch. Ration. Mech. Anal., 178 (2005),
  pp.~301--325.

\bibitem{BS2011}
{\sc L.~C. Berselli, P.~Miloro, A.~Menciassi, and E.~Sinibaldi}, {\em Exact
  solution to the inverse {W}omersley problem for pulsatile flows in
  cylindrical vessels, with application to magnetic particle targeting}, Appl.
  Math. Comput. (doi: 10.1016/j.amc.2012.11.071),  (2012, in press).

\bibitem{BR2010}
{\sc L.~C. Berselli and M.~Romito}, {\em On {L}eray's problem for almost
  periodic flows}, J. Math. Sci. Univ. Tokyo, 19 (2012), pp.~69--130.

\bibitem{B2009}
{\sc C.~D. Bertram}, {\em A numerical investigation of waves propagating in the
  spinal cord and subarachnoid space in the presence of a syrinx}, J. Fluids
  Struct., 25 (2009), pp.~1189--1205.

\bibitem{Canic2006}
{\sc S.~{\v C}ani{\'c}, J.~Tambaca, G.~Guidoboni, A.~Mikelic, C.~Hartley, and
  D.~Rosenstrauch}, {\em Modeling viscoelastic behavior of arterial walls and
  their interaction with pulsatile blood flow}, SIAM J. Appl. Math., 67 (2006),
  pp.~164--193.

\bibitem{Caro2009}
{\sc C.~G. Caro}, {\em Discovery of the role of wall shear in atherosclerosis},
  Arterioscler Thromb Vasc Biol., 29 (2009), pp.~158--161.

\bibitem{Carpinlioglu2001}
{\sc M.~{\"O}. {\c C}arpinlioglu and M.~Y. G{\"u}ndogdu}, {\em A critical
  review on pulsatile pipe flow studies directing towards future research
  topics}, Flow Meas. Instrum., 12 (2001), pp.~163--174.

\bibitem{CMZ1994}
{\sc G.~Chen, P.~J. Morris, and J.~Zhou}, {\em Visualization of special
  eigenmode shapes of a vibrating elliptical membrane}, SIAM Rev., 36 (1994),
  pp.~453--469.

\bibitem{cheng2012the}
{\sc S.~Cheng, M.~A. Stoodley, J.~Wong, S.~Hemley, D.~F. Fletcher, and L.~E.
  Bilston}, {\em {The presence of arachnoiditis affects the characteristics of
  CSF flow in the spinal subarachnoid space: a modelling study}}, J. Biomech.,
  45 (2012), pp.~1186--91.

\bibitem{CK2012}
{\sc S.~Cirovic and M.~Kim}, {\em A one-dimensional model of the spinal
  cerebrospinal-fluid compartment}, J. Biomech. Eng., 134 (2012), p.~021005.

\bibitem{dirocco2011hyd}
{\sc C.~Di~Rocco, P.~Frassanito, L.~Massimi, and S.~Peraio}, {\em
  {Hydrocephalus and Chiari type I malformation}}, Childs Nerv. Syst., 27
  (2011), pp.~1653--64.

\bibitem{E2012}
{\sc N.~S.~J. Elliott}, {\em Syrinx fluid transport: Modeling
  pressure-wave-induced flux across the spinal pial membrane}, J. Biomech.
  Eng., 134 (2012), p.~31006.

\bibitem{FQV2006}
{\sc L.~Formaggia, A.~Quarteroni, and A.~Veneziani}, {\em The circulatory
  system: from case studies to mathematical modeling}, in Complex systems in
  biomedicine, Springer Italia, Milan, 2006, pp.~243--287.

\bibitem{FVV2010}
{\sc L.~Formaggia, A.~Veneziani, and C.~Vergara}, {\em Flow rate boundary
  problems for an incompressible fluid in deformable domains: formulations and
  solution methods}, Comput. Methods Appl. Mech. Engrg., 199 (2010),
  pp.~677--688.

\bibitem{Abbott2010}
{\sc T.~W.~R. Fountain, P.~V. Kailat, and J.~J. Abbott}, {\em Wireless control
  of magnetic helical microrobots using a rotating-permanent-magnet
  manipulator}, in Proceedings of 2010 IEEE International Conference on
  Robotics and Automation, Anchorage, AK, USA, May 2010, pp.~576--581.

\bibitem{Gal1994a}
{\sc G.~P. Galdi}, {\em An introduction to the mathematical theory of the
  {N}avier-{S}tokes equations. {V}ol. {I}}, vol.~38 of Springer Tracts in
  Natural Philosophy, Springer-Verlag, New York, 1994.
\newblock Linearized steady problems.

\bibitem{GR2005}
{\sc G.~P. Galdi and A.~M. Robertson}, {\em The relation between flow rate and
  axial pressure gradient for time-periodic {P}oiseuille flow in a pipe}, J.
  Math. Fluid Mech., 7 (2005), pp.~S215--S223.

\bibitem{GBBCFMS2011}
{\sc A.~Gizzi, M.~Bernaschi, D.~Bini, C.~Cherubini, S.~Filippi, S.~Melchionna,
  and S.~Succi}, {\em Three-band decomposition analysis of wall shear stress in
  pulsatile flows}, Phys. Rev. E, 83 (2011), p.~031902.

\bibitem{GL2004}
{\sc L.~Greengard and J.-Y. Lee}, {\em Accelerating the nonuniform fast
  {F}ourier transform}, SIAM Rev., 46 (2004), pp.~443--454.

\bibitem{GPK2008}
{\sc S.~Gupta, D.~Poulikakos, and V.~Kurtcuoglu}, {\em {Analytical solution for
  pulsatile viscous flow in a straight elliptic annulus and application to the
  motion of the cerebrospinal fluid}}, {Phys. Fluids}, {20} ({2008}),
  p.~093607.

\bibitem{Haddad2010}
{\sc K.~Haddad, O.~Ertun\ifmmode~\mbox{\c{c}}\else \c{c}\fi{}, M.~Mishra, and
  A.~Delgado}, {\em Pulsating laminar fully developed channel and pipe flows},
  Phys. Rev. E, 81 (2010), p.~016303.

\bibitem{HZ1998}
{\sc M.~Haslam and M.~Zamir}, {\em {Pulsatile ﬂow in tubes of elliptic cross
  sections}}, Ann. Biomed. Eng., 26 (1998), p.~780.

\bibitem{Haverkort2009}
{\sc J.~W. Haverkort, S.~Kenjere, and C.~R. Kleijn}, {\em Magnetic particle
  motion in a {P}oiseuille flow}, Phys. Rev. E, 80 (2009), p.~016302.

\bibitem{Hentschel2010}
{\sc S.~Hentschel, K.-A. Mardal, A.~E. Lovgren, S.~Linge, and V.~Haughton},
  {\em {Characterization of Cyclic CSF Flow in the Foramen Magnum and Upper
  Cervical Spinal Canal with MR Flow Imaging and Computational Fluid
  Dynamics}}, Am. J. Neuroradiol., 31 (2010), pp.~997--1002.

\bibitem{Hoi2010}
{\sc Y.~Hoi, B.~A. Wasserman, Y.~Y.~J. Xie, S.~S. Najjar, L.~Ferruci, E.~G.
  Lakatta, G.~Gerstenblith, and D.~A. Steinman}, {\em Characterization of
  volumetric flow rate waveforms at the carotid bifurcations of older adults},
  Physiol. Meas., 31 (2010), pp.~291--302.

\bibitem{Huang2006}
{\sc H.-F. Huang and C.-L. Lai}, {\em Enhancement of mass transport and
  separation of species by oscillatory electroosmotic flows}, Proc. R. Soc.
  Lond. Ser. A Math. Phys. Eng. Sci., 462 (2006), pp.~2017--2038.

\bibitem{Irani2008}
{\sc D.~N. Irani}, {\em Cerebrospinal Fluid in Clinical Practice}, Saunders,
  2008.

\bibitem{Kha1957}
{\sc S.~R. Khamrui}, {\em On the flow of a viscous liquid through a tube of
  elliptic section under the influence of a periodic pressure gradient}, Bull.
  Calcutta Math. Soc., 49 (1957), pp.~57--60.

\bibitem{Korichi2009}
{\sc A.~Korichi, L.~Oufer, and G.~Polidori}, {\em Heat transfer enhancement in
  self-sustained oscillatory flow in a grooved channel with oblique plates},
  Int. J. Heat Mass Transfer, 52 (2009), pp.~1138--1148.

\bibitem{Krejza2006}
{\sc J.~Krejza, M.~Arkuszewski, S.~E. Kasner, J.~Weigele, A.~Ustymowicz, R.~W.
  Hurst, B.~L. Cucchiara, and S.~R. Messe}, {\em {C}arotid {A}rtery {D}iameter
  in {M}en and {W}omen and the {R}elation to {B}ody and {N}eck {S}ize}, Stroke,
  37 (2006), pp.~1103--5.

\bibitem{Lai2004}
{\sc M.-C. Lai}, {\em Fast direct solver for {P}oisson equation in a 2{D}
  elliptical domain}, Numer. Methods Partial Differential Equations, 20 (2004),
  pp.~72--81.

\bibitem{Linn2009}
{\sc A.~A. Linninger, M.~Xenos, B.~Sweetman, S.~Ponkshe, X.~Guo, and R.~Penn},
  {\em A mathematical model of blood, cerebrospinal fluid and brain dynamics},
  J. Math. Biol., 59 (2009), pp.~729--759.

\bibitem{LYA2001}
{\sc F.~Loth, M.~A. Yardimci, and N.~Alperin}, {\em Hydrodynamic modeling of
  cerebrospinal fluid motion within the spinal cavity}, J. Biomech. Eng., 123
  (2001), pp.~71--79.

\bibitem{loukas2011}
{\sc M.~Loukas, B.~J. Shayota, K.~Oelhafen, J.~H. Miller, J.~J. Chern, R.~S.
  Tubbs, and W.~J. Oakes}, {\em {Associated disorders of Chiari Type I
  malformations: a review}}, Neurosurg Focus, 31 (2011), p.~E3.

\bibitem{Abbott2011}
{\sc A.~W. Mahoney, J.~C. Sarrazin, E.~Bamberg, and J.~J. Abbott}, {\em
  Velocity control with gravity compensation for magnetic helical
  microswimmers}, Advanced Robotics, 25 (2011), pp.~1007--1028.

\bibitem{McLach1947}
{\sc N.~W. McLachlan}, {\em {Theory and application of Mathieu functions}},
  Clarendon Press, 1947.

\bibitem{Nelson2010}
{\sc B.~J. Nelson, I.~K. Kaliakatsos, and J.~J. Abbott}, {\em Microrobots for
  minimally invasive medicine}, Annu. Rev. Biomed. Eng., 12 (2010), pp.~55--85.

\bibitem{Peacock1998}
{\sc J.~Peacock, T.~Jones, C.~Tock, and R.~Lutz}, {\em The onset of turbulence
  in physiological pulsatile flow in a straight tube}, Exp. Fluids, 24 (1998),
  pp.~1--9.

\bibitem{Pil2007}
{\sc K.~Pileckas}, {\em The {N}avier-{S}tokes system in domains with
  cylindrical outlets to infinity. {L}eray's problem}, in Handbook of
  mathematical fluid dynamics, {V}ol. {IV}, North-Holland, Amsterdam, 2007,
  pp.~445--647.

\bibitem{Pop2000}
{\sc S.~B. Pope}, {\em Turbulent flows}, Cambridge University Press, Cambridge,
  2000.

\bibitem{Sarpeshkar2012}
{\sc B.~I. Rapoport, J.~T. Kedzierski, and R.~Sarpeshkar}, {\em A glucose fuel
  cell for implantable brainmachine interfaces}, PLoS ONE, 7 (2012),
  p.~e38436.

\bibitem{kumar2010}
{\sc Y.~V.~K. Ravi~Kumar, P.~S.~V.~H.~N. Krishna~Kumari, M.~V. Ramana~Murthy,
  and S.~Sreenadh}, {\em Unsteady peristaltic pumping in a finite length tube
  with permeable wall}, J. Fluids Eng., 132 (2010), p.~101201.

\bibitem{RD2004}
{\sc S.~Ray and F.~Durst}, {\em Semianalytical solutions of laminar fully
  developed pulsating flows through ducts of arbitrary cross sections}, Phys.
  Fluids, 16 (2004), p.~4371.

\bibitem{RSO2009}
{\sc A.~M. Robertson, A.~Sequeira, and R.~G. Owens}, {\em Rheological models
  for blood}, in Cardiovascular mathematics, vol.~1 of MS\&A. Model. Simul.
  Appl., Springer Italia, Milan, 2009, pp.~211--241.

\bibitem{RKHM2001}
{\sc M.~B. Robertson, U.~K{\"o}hler, P.~R. Hoskins, and I.~Marshall}, {\em Flow
  in elliptical vessels calculated for a physiological waveform}, J. Vasc.
  Res., 38 (2001), pp.~73--82.

\bibitem{rogers2010blood}
{\sc K.~Rogers}, {\em Blood: Physiology and Circulation}, Human Body, Rosen
  Publishing Group, 2010.

\bibitem{Sex1930}
{\sc T.~Sexl}, {\em {U}ber den von {E.~G.~R}ichardson entdeckten
  ``{A}nnulareffekt"}, Z. Phys, 61 (1930), pp.~179--221.

\bibitem{Shaffer2011}
{\sc N.~Shaffer, B.~Martin, and F.~Loth}, {\em Cerebrospinal fluid
  hydrodynamics in type {I} {C}hiari malformation}, Neurol. Res., 33 (2011),
  pp.~247--260.

\bibitem{SW2009}
{\sc J.~Shen and L.-L. Wang}, {\em On spectral approximations in elliptical
  geometries using {M}athieu functions}, Math. Comp., 78 (2009), pp.~815--844.

\bibitem{Trip2012}
{\sc R.~Trip, D.~J. Kuik, J.~Westerweel, and C.~Poelma}, {\em An experimental
  study of transitional pulsatile pipe flow}, Phys. Fluids, 24 (2012),
  p.~014103.

\bibitem{Verma1960ellAnn}
{\sc P.~D. Verma}, {\em The pulsating viscous flow superposed on the steady
  laminar motion of incompressible fluid between two co-axial cylinders}, Proc.
  Indian Acad. Sci. Math. Sci., 26 (1960), pp.~447--458.

\bibitem{Verma1960ell}
\leavevmode\vrule height 2pt depth -1.6pt width 23pt, {\em The pulsating
  viscous flow superposed on the steady laminar motion of incompressible fluid
  in a tube of elliptic section}, Proc. Indian Acad. Sci. Math. Sci., 26
  (1960), pp.~282--297.

\bibitem{wagshul2011the}
{\sc M.~E. Wagshul, P.~K. Eide, and J.~R. Madsen}, {\em The pulsating brain: A
  review of experimental and clinical studies of intracranial pulsatility},
  Fluids Barriers CNS, 8 (2011), p.~5.

\bibitem{Wang2005}
{\sc X.~Wang and N.~Zhang}, {\em {Numerical analysis of heat transfer in
  pulsating turbulent flow in a pipe}}, Int. J. Heat Mass Transfer, 48 (2005),
  pp.~3957--3970.

\bibitem{Wom1955}
{\sc J.~R. Womersley}, {\em {Method for the calculation of velocity, rate of
  flow and viscous drag in arteries when the pressure gradient is known}}, J.
  Physiol., 127 (1955), pp.~553--563.

\end{thebibliography}

\end{document}